\documentclass{amsart}
\usepackage[alphabetic]{amsrefs}
\usepackage{geometry}
\usepackage{amssymb}
\usepackage{mathtools}
\usepackage{bm}

\usepackage{csquotes}
\usepackage{tikz}
\usetikzlibrary{calc,positioning,decorations.markings}
\usepackage[dvipsnames]{xcolor}
\usepackage{verbatim}
\usepackage{hyperref}
\usepackage{tikz-cd}

\newcommand{\cO}{\mathcal O}

\newcommand{\NN}{\mathbb{N}}
\newcommand{\QQ}{\mathbb{Q}}
\newcommand{\ZZ}{\mathbb{Z}}
\newcommand{\RR}{\mathbb{R}}
\newcommand{\field}{k}
\newcommand{\Path}{\mathbf{Path}}
\newcommand{\Order}{\mathbf{K}}
\newcommand{\homotopy}{\mathbf{h}}
\newcommand{\graph}{\mathbf{G}}
\newcommand{\Dual}{\mathbf{Dual}}

\newcommand{\Mepsilon}{D^{BC}_{\epsilon}}
\newcommand{\FFH}{F^{FH}}

\DeclareMathOperator{\Pic}{Pic}
\DeclareMathOperator{\Cl}{Cl}
\DeclareMathOperator{\Spec}{Spec}
\DeclareMathOperator{\Hom}{Hom}
\DeclareMathOperator{\End}{End}
\DeclareMathOperator{\Tor}{Tor}
\DeclareMathOperator{\Cell}{Cell}
\DeclareMathOperator{\mathspan}{span}
\DeclareMathOperator{\codim}{codim}

\DeclareMathOperator{\pdim}{pdim}
\DeclareMathOperator{\img}{img}

\newcommand{\id}{\operatorname{id}}
\newcommand{\Quiv}{\operatorname{Quiv}}

\DeclareMathOperator{\hull}{hull}

\newtheorem{theorem}{Theorem}[section]
\newtheorem{proposition}[theorem]{Proposition}
\newtheorem{lemma}[theorem]{Lemma}
\newtheorem{corollary}[theorem]{Corollary}

\theoremstyle{remark}
\newtheorem{definition}[theorem]{Definition}

\newtheorem{example}[theorem]{Example}
\newtheorem{remark}[theorem]{Remark}

\newcommand{\set}[2]{\left\{#1\,:\,#2\right\}}

\title{Constructing toric diagonal modules with short linear resolutions}
\author{Jay Yang}
\address{Jay Yang: Department of Mathematics, Vanderbilt University,
  Nashville, Tennessee, 37240, United States of America;
  {\normalfont \texttt{jay.k.yang@vanderbilt.edu}}}


\begin{document}

\begin{abstract}
  For a complete simplicial toric variety $X$,
  we construct a collection of Cohen--Macaulay modules with linear resolutions
  where the associated sheaf is the diagonal.
  Through this work, we see that the resolutions constructed by Favero--Huang from homotopy path algebras
  and those constructed by Hanlon--Hicks--Lazarev can be viewed as two dual methods
  to construct resolutions of the diagonal starting from a compatible $\ZZ^n$-stratification.
\end{abstract}

\maketitle


\section{Introduction}

This paper presents a construction of a linear resolution for the diagonal over any complete simplicial toric variety $X$ of length $\dim X$. While a classic result of commutative algebra, first stated carefully by Eisenbud and Goto \cite{eisenbud-goto}, implies that a standard graded truncation will give a linear resolution when truncated past the regularity, the resulting modules will essentially always have resolutions that are far longer than $\dim X$. 

There are two major case where short linear resolutions of the diagonal have been given. First is the Picard rank 2 case given in \cite{brown-sayrafi}, a comparison to which is discussed in Section~\ref{sec:dual}. Second is the resolutions from \cite{FH25}, which \cite{FH25}*{Proposition~6.29} proves are linear subject to some constraints to prove that the module of interest is Koszul. Unlike previous results, we will need almost no constraints on the toric variety $X$ to give a linear resolution.
The other major collection of resolutions of the diagonal from recent work, those of \cite{HHL24}, are essentially never linear.

Fix a complete normal toric variety $X$ and let $n$ be the number of rays in the fan associated to $X$.
Let $S=\field[x_1,\ldots,x_n,y_1,\ldots,y_n]$ be the Cox ring of $X\times X$ and $I_\Delta\subseteq S$ the saturated toric ideal for the diagonal.
Let $M_X\subseteq \ZZ^n$ be the character lattice associated to $X$.
For any $\epsilon\in \RR^n$ we can define the $S$-module
\[\Mepsilon = S\cdot \set{\frac{x^{\lceil \bm{p}-\epsilon\rceil}}{y^{\lceil \bm{p}-\epsilon\rceil}}}{\bm{p}\in M_X}\subseteq
\field[x_1^{\pm},\ldots,x_n^{\pm},y_1^{\pm},\ldots,y_n^{\pm}]/I_\Delta.\]
The $BC$ in the name $D^{BC}_\epsilon$ refers to the barycentric subdivision which will be a key to the construction of a resolution of this module.
Then the main theorem is that this module is Cohen--Macaulay and for a generic $\epsilon$, the module has a linear resolution. Moreover, so long as the toric variety $X$ is simplicial, the associated sheaf to $\Mepsilon$ is the structure sheaf of the diagonal.

\begin{theorem}
  \label{thm:linear}
  If $X$ is a complete normal toric variety  and $\epsilon\in \RR^n$ then $\pdim \Mepsilon = \codim \Mepsilon = \dim X$ and for generic $\epsilon$, the module $\Mepsilon$ has a linear minimal free resolution.
  Moreover, if $X$ is simplicial and $\Delta\subseteq X\times X$ is the diagonal subvariety, then $\left(\Mepsilon\right)^{\sim} = \cO_{\Delta}$ for $\epsilon\in \RR_{\geq 0}^n$ satisfying $\epsilon_i<1$.
\end{theorem}

While we do not directly construct the minimal free resolution, $\Mepsilon$ itself will arise naturally as $H_0$ of a non-minimal resolution.
The approach is to start with an $\epsilon$-shifted version of the ceiling stratification for the character lattice $M_X$ as a sublattice of $\ZZ^{n}$ from \cite{stratification-paper}*{Definition~3.1}. This idea of using $\epsilon$-shifts was strongly inspired by \cite{And24}. However, the $S$-modules which are resolved are distinct from those in \cite{And24}.
From the $\epsilon$-shifted ceiling stratification, we construct a second stratification, now a compatible $\ZZ^{2n}$-stratification obtained by combining the Lawrence lift with a barycentric subdivision. This new stratification will give a highly nonminimal resolution of $\Mepsilon$. Then a slight coarsening of the cell structure will provide a simpler resolution of $\Mepsilon$ from which we can prove the existence of a linear resolution.

There has been significant recent work on linear resolutions in a multigraded context.
Of particular note is the connection between truncations, multigraded-regularity and quasi-linear
resolutions from \cite{bruce-cranton-heller-sayrafi}*{Theorem~A}.
But there has also been significant work on linear strands \cite{BE24} and specifically
$N_p$-style conditions in the case of weighted projective space \cite{BE25}.
Our discussion of linearity will center on the standard graded notion of linearity
or what in this multigraded setting \cite{BE24}*{Definition 1.2} defines as \emph{strongly linear}.
That is to say, a resolution is linear if the matrices for the differentials contain only linear forms.

As compared to previous resolutions of the diagonal, for each toric variety we get multiple modules as $\epsilon$ varies. In particular, restricting to $\epsilon\in \RR^n$ where the
hyperplanes $\{\bm{v} \in M_X\otimes_{\ZZ} \RR \,|\, \bm{v}_i - \epsilon_i = j \}$ for $i=1,\ldots,n$ and $j\in \ZZ$ intersect transversely partitions $\RR^n$ into open chambers. In each chamber, the module $\Mepsilon$ and the corresponding resolution are constant. Moreover, these chambers give a sufficient genericity condition for Theorem~\ref{thm:linear}. The chamber structure that results from this is a refinement as a polyhedral complex of the secondary fan of $X$.

While this chamber structure described above gives a criterion for linearity of the minimal free resolution, a stronger condition for linearity can be achieved by looking at a quiver. For an $\epsilon\in \RR^n$, associate a compatible $\ZZ^n$-stratification, denoted $\psi_\epsilon:M_X \otimes_{\ZZ} \RR\rightarrow \ZZ^n$. To this stratification, associate a quiver $Q_{\psi_\epsilon}$ with vertices indexed by the image of $\psi_{\epsilon}$ viewed as elements of $\ZZ^n/M_X\cong \Cl(X)$ and edges arising from paths through $M_X$.
The linearity of the resolution is controlled by the structure of the quiver $Q_{\psi_\epsilon}$. The precise statement of this is given by Proposition~\ref{prop:linear}.
A more precise definition of the quivers discussed here is presented in Definition~\ref{def:monomial-quiver}.

While the module $\Mepsilon$ is not a truncation in the classic sense, with an appropriate grading and degree shift, it can be presented as a submodule of $S/I_{\Delta}$. For this connection to truncations, coarsen the $\Cl(X\times X)\cong \Cl(X)\times \Cl(X)$ grading on $S$ and $S/I_\Delta$ to a $\Cl(X)$-grading via the map $\Cl(X)\times \Cl(X)\to \Cl(X)$ defined by $(\alpha_1,\alpha_2)\mapsto \alpha_1+\alpha_2$. Then in this grading, there exist submodules generated in sufficiently ample degrees with short linear resolutions.

\begin{corollary}
  \label{cor:truncation}
  For $\alpha\in \Cl(X)$ sufficiently ample, there exists a submodule of $S/I_\Delta$ generated in degree $\alpha$ in the $\Cl(X)$-grading with a linear resolution of length $\dim X$.
\end{corollary}

The next collection of results arises from a specific partial minimization of the resolutions of $\Mepsilon$, distinct from the minimization used to obtain the existence of linear resolutions. The minimization is chosen such that the resolution still arises from a stratification on a simplicial complex. The resulting complex has a direct combinatorial description with a structure reminiscent of the order complex.

\begin{definition}
  \label{def:dual-complex}
  Given a compatible $\ZZ^n$-stratification $\psi$, the \emph{dual complex of the stratification} is the simplicial complex $\Dual_\psi$ with vertices given by $\img \psi\subseteq \ZZ^n$ and where $\{\bm{a}_0,\ldots,\bm{a}_d\}\in \Dual_\psi$ if there exists a chain of cells $\sigma_{0}\prec \cdots \prec \sigma_d$ with $\psi(\sigma_k)=\bm{a}_k$. 
\end{definition}

\begin{theorem}
  \label{thm:dual-complex}
  The dual complex of $\psi_\epsilon$ is $M_X$-invariant and its quotient has the homotopy type of a torus.
  Moreover, there exists a compatible $\ZZ^{2n}$-stratification $\psi_\epsilon^{dual}$ on the dual complex such that the corresponding resolution of $S$-modules is a free resolution of $\Mepsilon$.
\end{theorem}

To define this compatible $\ZZ^{2n}$-stratification, a slight generalization of \cite{stratification-paper}*{Definition~2.6} is needed and is presented in Definition~\ref{def:stratification}.

Via this dual complex construction, we get a connection to the resolutions of \cite{FH25}.
Let $R$ be the Cox ring of $X$. Let $\FFH$ be the resolution of the diagonal bimodule given
by applying the constructions of \cite{FH25}*{Section~6}, to the algebra of $\Cl(X)$-graded endomorphisms
$A:=\End(T)_{\bm{0}}$ for the following $R$-module given by summing over the homomorphic image of
$\img \psi_\epsilon$ in $\ZZ^{n}/M_X\cong \Cl(X)$:
\[T := \bigoplus_{\alpha\in \img \psi_\epsilon / M_X} R(\alpha).\]
Then translating from $A$-bimodules to $R\otimes_{\field} R$-modules via a tensor against $T \boxtimes_{\field}\Hom(T,R)$ as an $A$-bimodule gives a complex of free modules:
\[\FFH\otimes_{A\otimes_{\field}A^{op}} (T\boxtimes_{\field} \Hom(T,R))\cong \Hom(T,R)\otimes_{A} \FFH \otimes_{A} T.\]

This basic construction is inspired by the tilting bundle arguments of the Coherent--Constructible Correspondence; see \cite{baer-tilting}*{Theorem 3.1.2} for comparison.
However $T$ is not usually a tilting bundle and so there is not an equivalence of categories between the appropriate derived categories.

Nonetheless, the resolution from the stratification on the dual complex does arise as a subcomplex of the above complex. To state the theorem, denote the resolution of $\Mepsilon$ that arises from the dual complex as $\pi_{\Lambda(M_X)}(F_{\psi_{\epsilon}^{dual}})$. This notation will be justified in Section~\ref{sec:background}. 

\begin{theorem}
  \label{thm:hpa-pdim}
  Fix $\epsilon \in \RR^n$.
  Let $T'=\Hom(T,R)$ and let $A=\End(T)_{\bm{0}}$.
  The complex $\FFH\otimes_{A\otimes_{\field}A^{op}} (T\boxtimes_{\field} T')$ contains the resolution $\pi_{\Lambda(M_X)}(F_{\psi_{\epsilon}^{dual}})$ as a subcomplex, and
  \[H_0(\FFH\otimes_{A\otimes_{\field}A^{op}} (T\boxtimes_{\field} T')) = H_0(\pi_{\Lambda(M_X)}(F_{\psi_{\epsilon}^{dual}})).\]
\end{theorem}

\section{Background}
\label{sec:background}

We begin with a short discussion of toric varieties and toric ideals and some associated notation. Throughout, $\field$ will be a field.
For a saturated sublattice $L\subseteq \ZZ^n$, define the toric ideal in $\field[x_1,\ldots,x_n]$ associated to $L$ by
\[I_L = \langle x^{\bm{v}_+}-x^{\bm{v}_{-}}\,|\, \bm{v}_{+},\bm{v}_{-}\in \NN^n\text{ and } \bm{v}_{+}-\bm{v}_{-}\in L\rangle.\]

A key observation is that $\field[x_1,\ldots,x_n]/I_L\cong \field [\NN^{n}/L]$ where $\NN^{n}/L$ is
the homomorphic image of $\NN^{n}$ under the quotient $\ZZ^n\rightarrow \ZZ^n/L$.
As the isomorphism is canonical, we will identify $\field[x_1,\ldots,x_n]/I_L$ with $\field [\NN^{n}/L]$.
As an abuse of notation, for the quotient by $I_L$ in the Laurent polynomial ring we will define
\[\field[x_1^{\pm},\ldots,x_n^{\pm}]/I_L:=\field[x_1^{\pm},\ldots,x_n^{\pm}]/(I_L\field[x_1^{\pm},\ldots,x_n^{\pm}]) \cong \field[\ZZ^n/L].\]
As with the quotient of the polynomial ring, we will identify $\field[x_1^{\pm},\ldots,x_n^{\pm}]/I_L$ with $\field[\ZZ^n/L]$.

Now let us fix some specific rings for the remainder of the paper. Throughout, $n$ will be a fixed integer and $S$ will be the polynomial ring $\field[x_1,\ldots,x_n,y_1,\ldots,y_n]$ with $R$ the polynomial ring $\field[x_1,\ldots,x_n]$. Observe in particular that $S\cong \field[\NN^{2n}]\cong R\otimes_{\field} R$ and $R\cong \field[\NN^{n}]$. Typically, $R$ will be the Cox ring of a toric variety $X$ and $S$ with be the Cox ring of $X\times X$ however, this is not required outside of examples until Section~\ref{sec:diagonal}.
In general, we will follow the conventions of \cite{CLS11} on toric varieties. In particular, the Cox ring of $X$ has a $\Cl(X)$-grading and a module over the Cox ring is homogeneous if it is homogeneous with respect to the $\Cl(X)$-grading.

Throughout, for $\bm{u}\in \ZZ^n$ define $\left|\bm{u}\right|:=\sum_{i=1}^{n} \left|\bm{u}_i\right|$. Additionally, for $\bm{u},\bm{v}\in \ZZ^n$ denote by $\bm{u}\leq \bm{v}$ the partial order given by componentwise comparison and $\bm{u}<\bm{v}$ when $\bm{u}\leq \bm{v}$ and $\bm{u}\neq \bm{v}$. Moreover, assuming that $L\subseteq \ZZ^n$ is a saturated sublattice satisfying $L\cap \NN^n = \{0\}$, there is an induced partial order on $\ZZ^n/L$ where $\alpha\leq \beta$ if there exist $\bm{u}\in \alpha$ and $\bm{v}\in \beta$ such that $\bm{u}\leq \bm{v}$.

\subsection*{Cellular Complexes}

For this paper, a topological space $\mathcal{X}$ is said to have a \emph{cell structure} if it is homeomorphic to a CW complex in the sense of \cite{Hat02}*{}. Moreover, a cell of this cell structure is the image of a cell of the CW complex under the homeomorphism. Unless otherwise stated, the cells will be open cells. In particular 0-cells in $\mathcal{X}$ are points corresponding to vertices of the CW complex and each $k$-cell for $k>0$ is a subset of $\mathcal{X}$ homeomorphic to an open disc arising from open cells of the corresponding CW complex. Moreover, all cell structures in this paper will either be \emph{regular} in that the CW complex associated to it is regular, or be quotients of a regular cell structure under an action preserving the cell structure. Observe that if a cell structure is regular, then one can define a closed cell as the closure of an open cell and these closed cells will be homeomorphic to closed discs. Throughout cell complex and CW complex will be used synonymously.

If $\sigma,\tau$ are cells in a regular cell structure, then $\tau$ is a \emph{face} of $\sigma$ if $\tau$ is contained in the boundary of $\sigma$ and $\tau$ is a \emph{facet} of $\sigma$ if further $\tau$ is codimension 1 in the closure $\overline{\sigma}$. Denote that $\tau$ is a face of $\sigma$ by $\tau\prec \sigma$. Define $\tau\preceq \sigma$ to mean $\tau\prec \sigma$ or $\tau=\sigma$.

\subsection*{Stratifications and Quivers}

We begin by recalling the definition of a compatible $\ZZ^n$-stratification from \cite{stratification-paper}*{Definition~2.6}. As mentioned in the intro, we will need a slight generalization for Section~\ref{sec:dual}. Throughout, we give $\ZZ^n$ a topology via the Alexandrov topology with open upper sets, that is to say a subset $U$ of $\ZZ^n$ is open when $\bm{u}\leq \bm{v}$ and $\bm{u}\in U$ implies $\bm{v}\in U$.

\begin{definition}[\cite{stratification-paper}*{Definition~2.6}]
  \label{def:stratification}
  Fix a sublattice $L\subseteq \ZZ^n$ and endow $\ZZ^n$ with the Alexandrov topology under the componentwise order.
  Fix a topological space $V$ with a free $L$-action denoted $\bm{p}+\bm{v}$ for $\bm{p}\in V$ and $\bm{v}\in L$.
  A \emph{compatible $\ZZ^n$-stratification} on $V$ is a pair of a continuous map $\psi: V \rightarrow \ZZ^n$
  and an $L$-invariant cell structure on $V$ such that $\psi$ is constant on each open cell in $V$
  and satisfies $\psi(\bm{p}+\bm{v})=\psi(\bm{p})+\bm{v}$ for $\bm{p}\in V$ and $\bm{v}\in L$.
  For each open cell $\sigma\subseteq V$, define $\psi(\sigma):=\psi(\bm{p})$ where $\bm{p}$ is any point in $\sigma$.
\end{definition}

As in \cite{stratification-paper}, for $L\subseteq \ZZ^n$, denote by $L_{\RR}$ the vector space $L\otimes_{\ZZ} \RR$ viewed as a topological space with the standard Euclidean topology. On any cell $\sigma$, the value $\psi(\sigma)$ will also be called the \emph{weight} associated to the cell $\sigma$.

\begin{remark}
  The key propositions of \cite{stratification-paper}*{Section~2} go through with this modification as they only use the fact that $L_{\RR}$ has a free $L$-action by addition of vectors. In particular, \cite{stratification-paper}*{Definition 2.10} goes through without significant issue except that the complex now has length $\dim V$ instead of $\dim L_{\RR}$. The key propositions of interest are \cite{stratification-paper}*{Proposition 2.13, Corollary 2.15, and Proposition 2.17} all of which only depend on the construction in \cite{stratification-paper}*{Definition 2.10} and that $L_{\RR}$ is a CW complex that is invariant under a free $L$-action.

  Moreover, except for the compatible $\ZZ^n$-stratification on the dual complex defined in Section~\ref{sec:dual}, we will use $V=L_{\RR}$ which exactly recovers \cite{stratification-paper}*{Definition~2.6}. With care, one could prove that the dual complex is in fact homeomorphic to $L_{\RR}$. As that is not necessary for the results of this paper, we avoid it.
\end{remark}

The main stratifications used in \cite{stratification-paper} are ceiling stratifications, for which the cell structure arises from an infinite $L$-invariant hyperplane arrangement.
The starting point for this paper is an $\epsilon$-shifted variant.
We will call these \emph{shifted ceiling stratifications}. 
For $L\subseteq \ZZ^n$ and $\epsilon \in \RR^n$, define $\psi_\epsilon(\bm{v}) = \lceil \bm{v} - \epsilon \rceil$.
The underlying cell structure can be constructed in a similar way to \cite{stratification-paper}*{Remark~2.1}, where the maximal cells are the chambers of the hyperplane arrangement with hyperplanes of the form
\[\mathcal{H}_{i,j}=\{\bm{p}\in L_{\RR} \,|\, \bm{p}_i-\epsilon_i=j\}\text{ for } i=1,\ldots,n \text{ and } j\in \ZZ.\]
The case where $\epsilon=\bm{0}$ exactly recovers the ceiling stratification case.

The next two lemmas represent a couple of structural statements about the ceiling stratification which will be important later.
\begin{lemma}
  \label{lem:strata-maximal-cells}
  If $\psi$ is a shifted ceiling stratification, then the subsets
  \[(L_{\RR})_{\bm{u}}=\{\bm{p}\in L_{\RR} \,|\, \psi(\bm{p}) = \bm{u}\},\]
  called \emph{strata}, each contain a unique cell maximal with respect to inclusion among cells in $(L_{\RR})_{\bm{u}}$.
\end{lemma}
\begin{proof}
  Suppose that $\sigma$ and $\tau$ are maximal cells of $(L_{\RR})_{\bm{u}}$. As both cells are maximal, neither may be contained in the closure of the other.
  As the cell structure on $L_{\RR}$ arises from the chambers of a hyperplane arrangement, there exists a hyperplane of the stratification separating $\sigma$ and $\tau$.
  But the existence of such a hyperplane implies that $\psi(\sigma)\neq \psi(\tau)$, so they cannot both be cells of $(L_{\RR})_{\bm{u}}$.
\end{proof}

In addition to the strata, for any compatible $\ZZ^n$-stratification, we also define two other subsets of $L_{\RR}$ associated to a degree $\bm{u}\in \ZZ^n$:
\[(L_{\RR})_{\leq \bm{u}}:=\{\bm{v}\in L_{\RR} \,|\, \psi(\bm{v})\leq \bm{u}\}
\text{ and }
(L_{\RR})_{< \bm{u}}:=\{\bm{v}\in L_{\RR} \,|\, \psi(\bm{v})< \bm{u}\}\]

\begin{lemma}
  \label{lem:unique-maximal-cell}
  Let $\psi$ be a shifted ceiling stratification and $\bm{v}\leq \bm{u}$. Let $\tau$ be a cell with $\psi(\tau)=\bm{u}$.
  If there exists a cell $\sigma\preceq \tau$ such that $\psi(\sigma)=\bm{v}$, then there exists a unique cell satisfying $\psi(\sigma')=\bm{v}$ and $\sigma'\preceq \tau$ maximal with respect to inclusion.
\end{lemma}
\begin{proof}
  In the case where $\bm{v}=\bm{u}$, $\tau\preceq\tau$ and $\psi(\tau)=\bm{v}$ so $\sigma'=\tau$ is the unique maximum.
  Otherwise let $\sigma\prec \tau$ be a cell such that $\psi(\sigma)=\bm{v}$.
  By Lemma~\ref{lem:strata-maximal-cells} there exists a unique cell $\tau'$ in $(L_{\RR})_{\bm{v}}$ which is maximal with respect to inclusion.

  The cell structure on $L_{\RR}$ is polyhedral and the intersection $\overline{\tau}\cap \overline{\tau'}$
  is non-empty as the cell $\sigma$ satisfies $\sigma\preceq \tau$ and $\sigma\preceq \tau'$.
  Consequently, there exists a unique maximal cell $\sigma'$ in the intersection $\overline{\tau}\cap \overline{\tau'}$.
  Moreover, as $\sigma\subseteq \overline{\tau}\cap\overline{\tau'}$, it follows that $\sigma\preceq \sigma'$.
  Since $\psi(\sigma)=\psi(\tau')=\bm{v}$ this implies $\psi(\sigma')=\bm{v}$.
  By construction $\sigma'$ is maximal among faces of $\tau$ satisfying $\psi(\sigma')=\bm{v}$.
\end{proof}
We now define a pair of related quivers arising from the stratification.

\begin{definition}
  \label{def:monomial-quiver}
  Fix a polynomial ring $R$. A \emph{quiver labeled by $R$} is a pair of a quiver and an assignment of a monomial in $R$ to each edge of the quiver.
  
  Fix a sublattice $L\subseteq \ZZ^n$ and a compatible $\ZZ^n$-stratification $\psi$ supported on a regular cell complex with an $L$-action.
  Define the \emph{monomial quiver associated to $\psi$}, denoted $\widetilde{Q}_\psi$, as the quiver labeled by $R$ with vertices $\img\psi\subseteq \ZZ^n$ where for $\bm{u},\bm{v}\in\img \psi$ there is a unique edge $\bm{u}\rightarrow \bm{v}$ in $\widetilde{Q}_\psi$ if $\bm{v} < \bm{u}$ in the componentwise order and there exists a pair of cells $\tau\prec \sigma$ with $\tau$ a facet of $\sigma$ such that $\psi(\sigma)=\bm{u}$ and $\psi(\tau)=\bm{v}$. The labels of the edges are given by $x^{\bm{u}-\bm{v}}\in \field[x_1,\ldots,x_n]$.

  Given the quiver $\widetilde{Q}_{\psi}$, define the \emph{toric quiver} $Q_\psi$ to be the quiver labeled by $R$ arising as the quotient of $\widetilde{Q}_\psi$ under the action induced by translation by $L$. In particular, the vertices are the corresponding elements of $\ZZ^n/L$ and the edges are equivalence classes of edges from $\widetilde{Q}_\psi$. 
\end{definition}

We also fix some notation for paths in quivers. For a path $p$ in a quiver from $a$ to $b$, denote by $h(p)=b$ and $t(p)=a$ the head and tail of $p$ respectively.

\begin{lemma}
  \label{lem:quiver-basics}
  Let $L\subseteq \ZZ^n$ be a lattice, let $\psi$ be a compatible $\ZZ^n$-stratification on $L_{\RR}$, and let $R=\field[x_1,\ldots,x_n]$. Then the following statements hold:

  \begin{enumerate}
  \item All paths between a pair of vertices in $\widetilde{Q}_\psi$ have the same product of the labels.
  \item The quiver $Q_\psi$ has a well-defined labeling by $R$.
  \item Furthermore, if $L \cap \NN^n=\{0\}$ then $Q_\psi$ has no directed cycles.
  \end{enumerate}
\end{lemma}

\begin{proof}
  For any path from $\bm{u}$ to $\bm{v}$, the product of the monomials must be $x^{\bm{u}-\bm{v}}$ and so is independent of path.
  
  The quiver $\widetilde{Q}_\psi$ acquires its $L$-equivariance from the $L$-equivariance of $\psi$. That is to say, for every element of $\bm{v}\in L$, there exists an isomorphism of quivers $\widetilde{Q}_\psi \xrightarrow{\cong} \widetilde{Q}_\psi$ given by adding $\bm{v}$ to every vertex of $\widetilde{Q}_\psi$. Moreover, the edges retain the same labels under this isomorphism. This $L$-equivariance implies that the quotient quiver $Q_\psi$ is well-defined and the edge labels in $Q_\psi$ are well-defined.

  Now assume that $L\cap \NN^n=\{0\}$. The monotonicity of the edges of the quiver implies that the quiver $\widetilde{Q}_\psi$ does not have any directed cycles.
  So to show that $Q_\psi$ is acyclic, it suffices to show that for any vertex $\bm{v}$ and any $L$-translate $\bm{v}'$, there does not exist a directed path from $\bm{v}$ to $\bm{v}'$.
  However, as $L\cap \NN^n=\{0\}$ the vertices $\bm{v}$ and $\bm{v}'$ must be incomparable in the componentwise order.
  Thus there does not exist a directed path from $\bm{v}$ to $\bm{v}'$ and so $Q_\psi$ has no directed cycles.
\end{proof}

\begin{proposition}
  \label{prop:degree-quiver}
  Fix a sublattice $L\subseteq \ZZ^n$ satisfying $L \cap \NN^{n}=\{0\}$.
  Fix a shifted ceiling stratification $\psi$ for $L$.
  If $\alpha,\beta$ are vertices of $Q_\psi$ such that $\beta<\alpha$ in the induced order on $\ZZ^n/L$ then there exists a path from $\alpha$ to $\beta$ in $Q_\psi$.
\end{proposition}
\begin{proof}
  If $\beta<\alpha$ in the induced order on $\ZZ^n/L$ then there exist $\bm{a},\bm{b}\in \ZZ^n$ satisfying $\bm{b}<\bm{a}$ with $\bm{a}\in \alpha$ and $\bm{b}\in\beta$. 
  Thus it suffices to show that for $\bm{a},\bm{b}\in \ZZ^n$ satisfying $\bm{b}<\bm{a}$, there exists a path in $\widetilde{Q}_\psi$ from $\bm{a}$ to $\bm{b}$.
  
  As any increasing chain in $\ZZ^n$ from $\bm{b}$ to $\bm{a}$ is of finite length, we can reduce to the case where $\bm{b},\bm{a}$ are vertices of $\widetilde{Q}_{\psi}$ such that there does not exist $\bm{u}\in \img \psi$ such that $\bm{b}< \bm{u} < \bm{a}$.
  
  Let $\bm{w}_{\bm{a}}$ and $\bm{w}_{\bm{b}}$ be points in $L_{\RR}$ such that $\psi(\bm{w}_{\bm{a}})=\bm{a}$ and $\psi(\bm{w}_{\bm{b}})=\bm{b}$.
  Then consider the line segment between $\bm{w}_{\bm{a}}$ and $\bm{w}_{\bm{b}}$.
  As $\psi$ is a shifted ceiling stratification, $\psi$ is monotonic on the line segment.
  Then $\psi$ is always equal to either $\bm{a}$ or $\bm{b}$ on the line segment as no $\bm{u}\in \img\psi$ satisfies $\bm{b}< \bm{u} < \bm{a}$.
  Since $\bm{b} < \bm{a}$, and $\psi$ is continuous with respect to the Alexandrov topology,
  there must exist a point $\bm{w}$ on the line segment such that $\psi(\bm{w})=\bm{b}$ with $\bm{w}\in \overline{(L_{\RR})_{\bm{a}}}$.
  Then let $\sigma$ be maximal among cells in $\overline{(L_{\RR})_{\bm{a}}}$ satisfying $\bm{w}\in \overline{\sigma}$ and $\psi(\sigma)=\bm{b}$.
  As $\bm{w}\in \overline{(L_{\RR})_{\bm{a}}}$, there must exist a cell $\sigma'$ such that $\psi(\sigma')=\bm{a}$ and $\sigma\prec \sigma'$.
  Choose the minimal such $\sigma'$. Then $\sigma$ is a facet of $\sigma'$.
  As such, there exists a pair of cells $\sigma\prec\sigma'$ with $\sigma$ a facet of $\sigma'$ satisfying $\psi(\sigma)=\bm{b}$ and $\psi(\sigma')=\bm{a}$.
  Thus there exists an edge $\bm{a}\rightarrow \bm{b}$ in the quiver $\widetilde{Q}_{\psi}$.
\end{proof}

We make two observations about the quiver $Q_\psi$ which are not critical to the proofs in this paper. First is that the quiver $Q_\psi$ is independent of the choice of cell structure on $\psi$. Second is that if $\psi$ is the ceiling stratification on $(M_X)_{\RR}$ for a toric variety $X$, the vertices for $Q_\psi$ are the $\Cl(X)$-graded degrees of the Bondal--Thomsen collection and the quiver can be viewed as arising from the endomorphisms of the Bondal--Thomsen collection.

\subsection*{Cellular Resolutions of Toric Ideals}

Let $V$ be a space with a free $L$-action. For a compatible $\ZZ^n$-stratification $\psi$ on $V$, associate a complex of modules $F_{\psi}$ as defined in \cite{stratification-paper}*{Definition~2.10}. Moreover, assuming that the stratification satisfies the conditions of \cite{stratification-paper}*{}, the complex $F_\psi$ resolves the module $M_\psi:=S\cdot \{x^{\psi(\bm{v})} \,|\, \bm{v} \in V\}$. However, the module $M_\psi$ is usually not the module of interest. In particular, there is an $L$-equivariant structure on $M_\psi$ and $F_\psi$ which one can quotient by. Using the notation of \cite{stratification-paper}, define the ring $S[L]$ by
\[S[L] := S\cdot \{z^{\bm{u}} \,|\, \bm{u}\in L\} \subseteq S[z_1^{\pm},\ldots,z_n^{\pm}].\]
Then the module $M_\psi$ has an $S[L]$-module structure arising from the $L$-action. Under the isomorphism $S\cong S[L]/\left<z^{\bm{a}}-1 \,|\, \bm{a}\in L\right>$, the $S$-module $M_\psi\otimes_{S[L]}S$ is resolved by the complex $F_\psi \otimes_{S[L]} S$. For brevity, denote by $\pi_L$ the functor $-\otimes_{S[L]} S$.

We also recall the direct construction of $\pi_L(F_\psi)$ discussed in \cite{stratification-paper}*{Remark~2.14}. For a stratification $\psi$, the complex $\pi_L(F_\psi)$ has summands given by
\[\pi_L(F_{\psi})_{i} := \bigoplus_{\dim \sigma = i}S(-\psi(\sigma)).\]
Fix a cell $\sigma$ in $V/L$ and fix $\sigma'$ a representative in $V$. Then differentials are given, up to sign, by
\[\pi_L(\partial_{\psi})_{i}(\sigma) := \sum_{\tau\text{ a facet of } \sigma}\left(\sum_{\tau'\in \tau\text{ and }\tau'\prec \sigma'} \pm x^{\psi(\sigma')-\psi(\tau')}\right)\tau.\]
For both $F_\psi$ and $\pi_L(F_\psi)$, denote by $S\sigma$ the summand of the complex corresponding to the cell $\sigma$.

\section{Subdividing Stratifications}
\label{sec:subdivision}

As remarked in the introduction, a key tool in this construction of short linear resolutions is the barycentric subdivision.
For this it will be convenient to restrict to a class of stratifications where the underlying cell complexes are polyhedral complexes. However, essentially all of the work here translates without issue so long as the cell structure on $L_{\RR}$ is regular.

\begin{definition}
  \label{def:polyhedral}
  A compatible $\ZZ^n$-stratification on $L_{\RR}$ is said to be \emph{polyhedral} if the underlying cell complex is a polyhedral complex.
\end{definition}

It is immediate that all shifted ceiling stratifications are polyhedral. We will present as much as possible of the subdivision work in this generality.

For a polyhedral cell complex, its \emph{barycentric subdivision} is a refinement of the cell complex. The vertices of the barycentric subdivision are given by the barycenters of each cell $\sigma\subseteq L_{\RR}$, which we will denote $v_{\sigma}$. The cells are given by the convex hulls of collections of vertices $\{v_{\sigma_0},\ldots,v_{\sigma_d}\}$ where 
\[\sigma_0\prec \sigma_1 \prec \cdots \prec \sigma_d.\]
In particular, the maximal cells of the barycentric subdivision are in bijection with maximal chains of cells from $L_{\RR}$.
Moreover, for any cell in the barycentric subdivision, the faces of that cell are in bijection with subchains.
For the remainder of the paper, when we write that a cell $\tau$ in the barycentric subdivision is given by $\hull\{v_{\sigma_0},\ldots,v_{\sigma_d}\}$
we assume that the cells are ordered such that $\sigma_0 \prec \cdots \prec \sigma_d$.

This identification of cells in the barycentric subdivision with chains of cells in the original complex
is equivalent to observing that the barycentric subdivision is a geometric realization of
the order complex of the face poset of our cell complex.
In fact this description gives an extension of the barycentric subdivision to regular cell complexes.

Using the same stratification function with this new cell structure is not productive, as doing so simply results in a chain complex that minimizes to the original one. Instead, from the compatible $\ZZ^n$-stratification $\psi$ on $L_{\RR}$, we will construct a compatible $\ZZ^{2n}$-stratification on $\Lambda(L)_{\RR}$ where \[\Lambda(L) := \{(\bm{v},-\bm{v}) \,|\, \bm{v}\in L\}\] is the Lawrence lift of $L$ as defined in \cite{SturmfelsGB}*{Equation~7.1}. In particular, the isomorphism of lattices $L\xrightarrow{\cong}\Lambda(L)$ by $\bm{v}\mapsto (\bm{v},-\bm{v})$ gives a homeomorphism $L_{\RR} \rightarrow \Lambda(L)_{\RR}$ and so we can use the barycentric subdivision as the cell structure on $\Lambda(L)_{\RR}$.

\begin{definition}
  Fix a polyhedral stratification $\psi$ on $L$.
  Give $\Lambda(L)_{\RR}$ the cell structure of the barycentric subdivision of $L_{\RR}$
  induced by the homeomorphism $L_{\RR}\rightarrow \Lambda(L)_{\RR}$.

  Denote by $\mathbf{max}$ the componentwise maximum of a collection of vectors.
  Then define the \emph{barycentric diagonal stratification} to be the compatible $\ZZ^{2n}$-stratification
  $\psi^{BC}$ on $\Lambda(L)_{\RR}$ where for a simplex $\tau=\hull\{v_{\sigma_0},\ldots,v_{\sigma_d}\}$ in the barycentric subdivision
  \[\psi^{BC}(\tau) := \mathbf{max}\left\{(\psi(\sigma_0),-\psi(\sigma_0)), \ldots, (\psi(\sigma_d),-\psi(\sigma_d))\right\}.\]
\end{definition}

Where it is unambiguous, we will assume that the cell structure on $\Lambda(L)_{\RR}$ is that acquired from the barycentric subdivision of $L_{\RR}$.
Otherwise, when needed for clarity we will denote by $\Lambda(L)_{\RR}^{BC}$ the space $\Lambda(L)_{\RR}$ with the cell structure arising as the barycentric subdivision of $L_{\RR}$.

\begin{lemma}
  If $\psi$ is a polyhedral compatible $\ZZ^{n}$-stratification on $L_{\RR}$ then $\psi^{BC}$ is a compatible $\ZZ^{2n}$-stratification of $\Lambda(L)$.
\end{lemma}
\begin{proof}
  It is immediate that $\psi^{BC}$ is constant on cells of the cell structure.
  The equivariance under the translation action of $\Lambda(L)$ follows from
  the equivariance of $\psi$ under the translation action of $L$.
  It remains to show that $\psi^{BC}$ is a continuous map when equipping $\ZZ^{2n}$ with the Alexandrov topology.
  For this it suffices to show that the set $(\Lambda(L)_{\RR})_{\leq \bm{u}}$ is closed for all $\bm{u}\in \ZZ^{2n}$.

  For $\tau=\hull\{v_{\sigma_0},\ldots,v_{\sigma_d}\}$, the value $\psi^{BC}(\tau)$ is the max of the values of $\psi^{BC}$ on all of the vertices in $\tau$. As any face of $\tau$ has a subset of the vertices of $\tau$ it follows that if $\tau'\preceq \tau$ then $\psi^{BC}(\tau')\leq \psi^{BC}(\tau)$. As such, the set $(\Lambda(L)_{\RR})_{\leq \bm{u}}$ is closed.
\end{proof}

As it will be convenient, we will also observe the following simpler equation for the value of $\psi^{BC}$ on cells of $\Lambda(L)_{\RR}$.
\begin{lemma}
  Let $\psi$ be a polyhedral stratification and let $\sigma_0\prec \cdots \prec \sigma_d$ be a chain of cells in $L_{\RR}$.
  For $\tau=\hull\{v_{\sigma_0},\ldots,v_{\sigma_d}\}\subseteq L_{\RR}^{BC}$,
  \[\psi^{BC}(\tau) = (\psi(\sigma_d),-\psi(\sigma_0)).\]
\end{lemma}
\begin{proof}
  Since $\sigma_0\prec \cdots \prec \sigma_d$, it follows that $\psi(\sigma_i)\leq \psi(\sigma_{i+1})$
  where $\ZZ^n$ is given the componentwise partial order.
  It follows that $\mathbf{max}(\psi(\sigma_0),\ldots,\psi(\sigma_d))=\psi(\sigma_d)$ and $\mathbf{max} (-\psi(\sigma_0), \ldots, -\psi(\sigma_d))=-\psi(\sigma_0)$. Thus
  \[\psi^{BC}(\tau) = \mathbf{max}\left\{(\psi(\sigma_0),-\psi(\sigma_0)), \ldots, (\psi(\sigma_d),-\psi(\sigma_d))\right\} = (\psi(\sigma_d),-\psi(\sigma_0)).\]
\end{proof}

\subsection*{Complexes from Subdivision}

Recall that $S=\field[x_1,\ldots,x_n,y_1,\ldots,y_n]\cong \field[\NN^{2n}]$.
Then the compatible $\ZZ^{2n}$-stratification $\psi^{BC}$ of $\Lambda(L)_{\RR}$ arising from the compatible $\ZZ^n$-stratification $\psi$ on $L_{\RR}$
gives a complex of $S$-modules $\pi_{\Lambda(L)}(F_{\psi^{BC}})$.
As an example, we can subdivide the ceiling stratification for the lattice $M_{\mathbb{P}^2}\subseteq \ZZ^{3}$.

\begin{example}
  Let $\psi$ be the ceiling stratification for the lattice $L=\mathspan_{\ZZ} \{(1,-1,0),(0,1,-1)\}\subseteq \ZZ^3$. Geometrically, this corresponds to a resolution of the structure sheaf of the torus identity point in $\mathbb{P}^2$. Observe that in this case $R=\field[x_1,x_2,x_3]$ has a $\ZZ$-grading via the isomorphism $\ZZ\cong \ZZ^3/L \cong \Cl(\mathbb{P}^2)$. From this grading, $S=\field[x_1,x_2,x_3,y_1,y_2,y_3]$ has a $\ZZ^2$-grading which we use below for our degree shifts.

  Subdividing the cell structure on $L_{\RR}$ yields a cell structure on $\Lambda(L)_{\RR}$ which in a fundamental domain has 6 vertices, 18 edges, and 12 2-cells, as illustrated in Figure~\ref{fig:bc-p2}.

  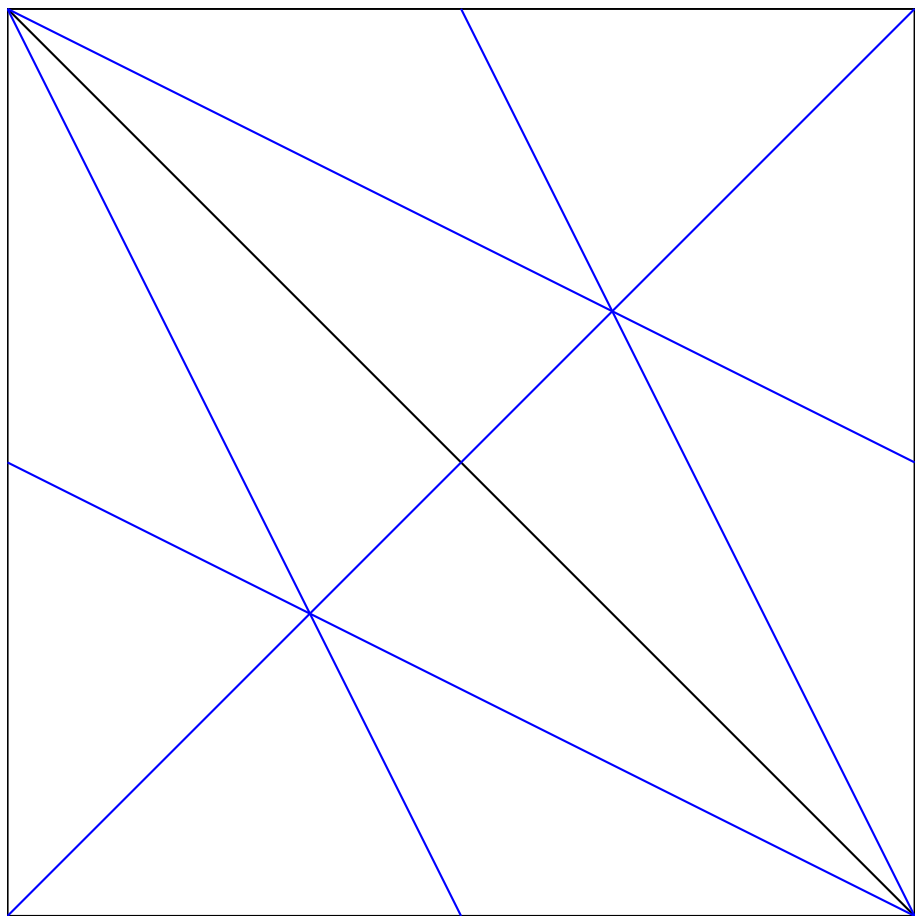
\begin{figure}
    \begin{tikzpicture}
      \clip (0,0) rectangle (6,6);
      
      \coordinate (origin) at (0,0);
      \coordinate (br) at (6,0);
      \coordinate (tr) at (6,6);
      \coordinate (tl) at (0,6);

      \draw ($(origin)-(0,1)$) -- ($(tl)+(0,1)$);
      \draw ($(origin)-(1,0)$) -- ($(br)+(1,0)$);
      \draw ($(br)-(0,1)$) -- ($(tr)+(0,1)$);
      \draw ($(tl)-(1,0)$) -- ($(tr)+(1,0)$);
      \draw ($(tl)-(1,-1)$) -- ($(br)+(1,-1)$);
      \draw ($(origin)-(2,-2)$) -- ($(origin)+(2,-2)$);
      \draw ($(tr)-(1,-1)$) -- ($(tr)+(1,-1)$);

      \coordinate (bm) at (3,0);
      \coordinate (lm) at (0,3);
      \coordinate (tm) at (3,6);
      \coordinate (rm) at (6,3);
      \coordinate (center) at (3,3);
      \coordinate (bc1) at (2,2);
      \coordinate (bc2) at (4,4);

      \begin{scope}[draw=blue]
        \draw (origin) -- (bc1);
        \draw (tl) -- (bc1);
        \draw (br) -- (bc1);
        \draw (bm) -- (bc1);
        \draw (lm) -- (bc1);
        \draw (center) -- (bc1);
        \draw (center) -- (bc2);
        \draw (tl) -- (bc2);
        \draw (br) -- (bc2);
        \draw (tr) -- (bc2);
        \draw (rm) -- (bc2);
        \draw (tm) -- (bc2);
      \end{scope}

    \end{tikzpicture}
    \caption{One fundamental domain of the barycentric subdivision of the cell complex associated to the ceiling stratification of $M_{\mathbb{P}^2}\subseteq \ZZ^{3}$.}
    \label{fig:bc-p2}
    
  \end{figure}

  This gives the following complex for $\pi_{\Lambda(L)}(F_{\psi^{BC}})$:

  \[
  S^1\oplus S(-1,1)^4 \oplus S(-2,2)^1 \xleftarrow{\partial_1}
  S(-1,0)^9\oplus S(-1,1)^3 \oplus S(-2,1)^3 \oplus S(-2,0)^3 \xleftarrow{\partial_2}
  S(-1,0)^6\oplus S(-2,0)^6
  \]

  With differentials given by
  \begin{align*}
    \partial_1 &= 
    \begin{psmallmatrix}
      x_1 & x_1 & x_2 & x_2 & x_3 & x_3 & x_1 & x_2 & x_3 & 0 & 0 & 0 & 0   & 0   & 0   & x_1x_2 & x_1x_3 & x_2x_3 \\
      -y_1 & 0 & -y_2 & 0   & 0  & 0   & 0 &  0  & 0  & 1 & 0 & 0  &0 & 0   &     x_3 &0&0&0\\
        0 & 0   & 0   & -y_2 & -y_3 & 0   & 0   & 0 & 0 & 0 & 1 & 0   & x_1    &0& 0   &0&0&0\\
        0 & -y_1 & 0   & 0   & 0  & -y_3 & 0   & 0   & 0 & 0& 0 & 1 & 0   &   x_2  &0&0&0&0\\
        0 & 0   & 0   & 0   & 0   & 0  & -y_1 & -y_2 & -y_3  & -1 & -1 & -1  & 0   & 0   & 0   &0&0&0\\
         0& 0   & 0   & 0   & 0   & 0  & 0 & 0 & 0 & 0   & 0   & 0   & -y_1 & -y_2 & -y_3 & -y_1y_2 & -y_1y_3 & -y_2y_3
    \end{psmallmatrix}\\
    \partial_2 &= 
    \begin{psmallmatrix}
      1   & 0   & 0   & 0   & 0   & 0   & 0   & 0   & 0   & x_3 & 0   & 0\\
      0   & 1   & 0   & 0   & 0   & 0   & 0   & x_2 & 0   & 0   & 0   & 0\\
      0   & 0   & 1   & 0   & 0   & 0   & 0   & 0   & 0   & 0   & 0   & x_3 \\
      0   & 0   & 0   & 1   & 0   & 0   & x_1 & 0   & 0   & 0   & 0   & 0\\
      0   & 0   & 0   & 0   & 1   & 0   & 0   & 0   & x_1 & 0   & 0   & 0\\
      0   & 0   & 0   & 0   & 0   & 1   & 0   & 0   & 0   & 0   & x_2 & 0\\
      -1  & -1  & 0   & 0   & 0   & 0   &0&0&0&0&0&0\\
      0   & 0   & -1  & -1  & 0   & 0   &0&0&0&0&0&0\\
      0   & 0   & 0   & 0   & -1  & -1  &0&0&0&0&0&0\\
      y_1 & 0   & y_2 & 0   & 0   & 0   &0&0&0&0&0&0\\
      0   & 0   & 0   & y_2 & y_3 & 0   &0&0&0&0&0&0\\
      0   & y_1 & 0   & 0   & 0   & y_3 &0&0&0&0&0&0\\
      0   & 0   & 0   & 0   & 0   & 0   & y_2 & 0 & y_3 &0&0&0\\
      0   & 0   & 0   & 0   & 0   & 0   & 0 &  y_1 & 0  & 0 & y_3 &0\\
      0   & 0   & 0   & 0   & 0   & 0   & 0 &  0   & 0  & y_1  & 0 & y_2 \\
      0   & 0   & 0   & 0   & 0   & 0   & -1& -1& 0 & 0 &0 &0\\
      0   & 0   & 0   & 0   & 0   & 0   & 0 & 0 & -1& -1& 0& 0\\
      0   & 0   & 0   & 0   & 0   & 0   & 0 & 0 & 0 &0 & -1& -1
    \end{psmallmatrix}
  \end{align*}

  One can check that the complex is a resolution of the module
  \[S\cdot \left\{1,\frac{x_1}{y_1},\frac{x_2}{y_2},\frac{x_1x_2}{y_1y_2},\frac{x_1x_2y_3}{x_3y_1y_2}\right\}
  = S\cdot \left\{1,\frac{x_1}{y_1},\frac{x_1x_2}{y_1y_2}\right\} \subseteq \field[\ZZ^6/\Lambda(L)].\]
\end{example}

The free resolution from the example can be minimized to the following:
\[S^1\oplus S(-1,1) \oplus S(-2,2) \leftarrow S(-1,0)^3 \oplus S(-2,1)^3 \leftarrow S(-2,0)^3.\]
One may observe that in the standard grading, this minimized resolution is linear.
The reasons for this are discussed more in Section~\ref{sec:diagonal}, where we provide a criterion for the linearity of this resolution.

We see in this example that the complex gives a resolution. To show this in general, by \cite{stratification-paper}*{Proposition~2.13}, we must show that the subcomplexes $(\Lambda(L)_{\RR})_{\leq \bm{u}}$ are acyclic for all $\bm{u}\in \ZZ^{2n}$.

\begin{proposition}
  \label{prop:resolution}
  If $\psi$ is a shifted ceiling stratification on $L_{\RR}$ for $L\subseteq \ZZ^n$, then the complex $\pi_{\Lambda(L)}(F_{\psi^{BC}})$ is a resolution of the $S$-module
  \[\pi_{\Lambda(L)}(M_{\psi^{BC}}) = S\cdot \{x^{\bm{a}}/y^{\bm{a}} \,|\, \bm{a}\in \img(\psi) \} \subseteq \field [\ZZ^{2n}/\Lambda(L)].\]
\end{proposition}

\begin{proof}[Proof of Proposition~\ref{prop:resolution}]
  Fix a lattice point $\bm{u}=(\bm{u}_1,\bm{u}_2)\in \ZZ^{2n}$ and define the subspace of $L_{\RR}$ given by
  \[\mathcal{U} := \{\bm{p} \in L_{\RR} \,|\, \psi(\bm{p})\leq \bm{u}_1\text{ and }-\psi(\bm{p})\leq \bm{u}_2\}.\]
  The isomorphism $L\cong \Lambda(L)$ gives a homeomorphism of topological spaces $\lambda:L_{\RR}\rightarrow \Lambda(L)_{\RR}$.

  The stratification $\psi^{BC}$ on $\Lambda(L)_{\RR}^{BC}$ has labels given by $(\psi(\lambda^{-1}(\bm{v})),-\psi(\lambda^{-1}(\bm{v})))$ for $\bm{v}$ a vertex of $\Lambda(L)_{\RR}^{BC}$. Moreover, the labels on all other cells are given by the componentwise maximums of the vertex labels for all the vertices contained in the cell. This implies that for all points $\bm{p}$ in $\Lambda(L)_{\RR}^{BC}$,\[\psi^{BC}(\bm{p})\geq (\psi(\lambda^{-1}(\bm{p})),-\psi(\lambda^{-1}(\bm{p}))).\] Consequently, \[(\Lambda(L)_{\RR})_{\leq \bm{u}}\subseteq \lambda(\mathcal{U}).\]

  Moreover, if $\sigma$ is a closed simplex of $\Lambda(L)_{\RR}^{BC}$ contained in $\lambda(\mathcal{U})$ and $\bm{v}$ is a vertex of $\sigma$, then $\psi^{BC}(\bm{v})\leq \bm{u}$.
  Thus $\psi^{BC}(\sigma)\leq \bm{u}$ and consequently $\sigma\in (\Lambda(L)_{\RR})_{\leq \bm{u}}$.
  Since $(\Lambda(L)_{\RR})_{\leq \bm{u}}$ is closed, it is the induced subcomplex of $\Lambda(L)_{\RR}^{BC}$ obtained from the union of all closed cells in $\lambda(\mathcal{U})$.

  We observe that if $\tau$ is an open simplex of $\Lambda(L)_{\RR}^{BC}$ contained in $\lambda(\mathcal{U})$ then at least one vertex of $\tau$ is in $\lambda(\mathcal{U})$. Then as $(\Lambda(L)_{\RR}^{BC})_{\leq \bm{u}}$ contains every closed cell, and thus every vertex in $\lambda(\mathcal{U})$, the open star of $(\Lambda(L)_{\RR}^{BC})_{\leq \bm{u}}$ inside of $\lambda(\mathcal{U})$ is exactly $\lambda(\mathcal{U})$. Observe that $\lambda(\mathcal{U})$ is contained in a finite subcomplex of $\Lambda(L)_{\RR}^{BC}$ and any induced subcomplex of a finite simplicial complex is a deformation retract of its open star \cite{Munkres}*{Lemma~70.1}. Thus $(\Lambda(L)_{\RR}^{BC})_{\leq \bm{u}}$ is a deformation retract of $\lambda(\mathcal{U})$. Finally, as $\mathcal{U}$ is polyhedral, it is contractible and so $(\Lambda(L)_{\RR}^{BC})_{\leq \bm{u}}$ is contractible and thus acyclic.
\end{proof}

\section{Diagonal Resolutions}
\label{sec:diagonal}

Now that we have a resolution arising from the barycentric subdivision, we will move on to its connection to the diagonal.
We recall that throughout $X$ is a complete normal toric variety with Cox ring $R=\field[x_1,\ldots,x_n]$ and $S=\field[x_1,\ldots,x_n,y_1,\ldots,y_n]\cong R\otimes_{\field} R$ is the Cox ring for $X\times X$.
Since the primary case of interest is the diagonal, for the remainder of the paper we restrict to the case where $L=M_X\subseteq \ZZ^n$. In this setting, the binomial ideal $I_{\Lambda(L)} = I_{\Delta}$ and so cuts out the diagonal subvariety $\Delta\subseteq X\times X$.

Our first result will be a simple condition to detect whether a homogeneous module gives the diagonal. Inside of $S$, define the ideal
\[I:= \left<x_i-y_i \,|\, i=1,\ldots,n \right>.\]
This ideal is not homogeneous with respect to the $\Cl(X\times X)$-grading on $S$.
However, $V(I)\subseteq \Spec(S)$ under the quotient map $\Spec(S)\dashrightarrow X\times X$ has image $\Delta$, specifically it cuts out the diagonal of the map on affine cones. Then to check whether two modules supported on the diagonal agree as sheaves, it suffices to check on $V(I)\subseteq \Spec(S)$.
\begin{lemma}
  \label{lem:homogenize-diagonal}
  Let $R=\field[x_1,\ldots,x_n]$ be the Cox ring of a simplicial toric variety $X$ and let $S=\field[x_1,\ldots,x_n,y_1,\ldots,y_n]$ be the Cox ring of $X\times X$. Let $\Delta\subseteq X\times X$ be the diagonal subvariety. Fix \[I:=\langle x^{\bm{u}}-y^{\bm{u}} \,|\, \bm{u}\in \mathbb{N}^n\rangle.\]
  
  Let $M$ be a homogeneous module over $S$ generated by elements in $\Pic(X\times X)$ degrees. If $\operatorname{ann}(M)=I_\Delta$ and $M/IM \cong S/I$ then $M^{\sim}\cong \cO_{\Delta}$.
\end{lemma}
\begin{proof}
  Since $\Cl(X\times X)\cong \Cl(X)\times \Cl(X)$, we write all degrees in $\Cl(X\times X)$ as a pair of degrees in $\Cl(X)$.
  In particular, where a single $\bm{0}$ is used as a subscript, it is with respect to the $\Cl(X)$-grading. Recall also that every $\Cl(X\times X)$-graded $S$-module has a $\Cl(X)$-graded structure from coarsening the grading via the map $\Cl(X)\times \Cl(X)\rightarrow \Cl(X)$ given by $(\alpha,\beta)\mapsto \alpha+\beta$. Let $\Sigma_{X\times X}$ and $\Sigma_{X}$ be fans of $X\times X$ and $X$ respectively

  We will check that the two sheaves agree using the standard affine cover of a toric variety associated to the fan for a toric variety.
  Since $\operatorname{ann} M = I_{\Delta}$, it suffices to check on
  affine charts corresponding to cones in $\Sigma_{X\times X}$ of the form
  $\tau = \tau'\times \tau'$ where $\tau'$ is a cone in $\Sigma_X$.
  Denote by $U_\tau$ and $U_{\tau'}$ the corresponding affine opens of $X\times X$ and $X$ respectively.
  Denote by $x_{\widehat{\tau}}$ and similarly $x_{\widehat{\tau'}}$
  the monomials given by the product of the variables corresponding to rays not in $\tau$ or $\tau'$.
  Then we explicitly compute the diagonal via the pushforward of the structure sheaf along the map
  $\varphi : X\hookrightarrow X\times X$.
  \begin{align*}
    \Gamma(U_\tau,\varphi_*(\cO_X)) &= \Gamma(U_{\tau'},\cO_X)\\
    &= \Gamma(U_{\tau'}, R^{\sim})\\
    &= R[x_{\widehat{\tau'}}^{-1}]_{\bm{0}} \\
    &\cong (M/IM)[x_{\widehat{\tau}}^{-1}]_{\bm{0}}\\
    &= \left(\bigoplus_{\alpha\in \Cl(X)}M[x_{\widehat{\tau}}^{-1}]_{(\alpha,-\alpha)}\right)/(IM[x_{\widehat{\tau}}^{-1}])_{\bm{0}}.
  \end{align*}

  As $X$ is simplicial, the subgroup of $\Cl(X\times X)$ generated the degrees of the variables
  appearing in $x_{\widehat{\tau}}$ contains $\Pic(X\times X)$.
  As $M$ is generated in $\Pic(X\times X)$ degrees, the module $M[x_{\widehat{\tau}}^{-1}]$
  is generated in degree $(\bm{0},\bm{0})$ as an $S[x_{\widehat{\tau}}^{-1}]$-module.
  Now for let $x^{u}y^{v}f\in M[x_{\widehat{\tau}}^{-1}]_{(\alpha,-\alpha)}$ where $f\in M[x_{\widehat{\tau}}^{-1}]_{(\bm{0},\bm{0})}$ and $x^{u}y^{v}\in S[x_{\widehat{\tau}}^{-1}]$.
  Then observe that $x^{u}-y^{u}\in I[x_{\widehat{\tau}}^{-1}]$ and so $x^{u}y^{v}f-y^{u+v}m\in IM[x_{\widehat{\tau}}^{-1}]$ and $\deg(y^{u+v}f) = (\bm{0},\bm{0})$.
  This implies that the quotient is spanned by elements of
  $M[x_{\widehat{\tau}}^{-1}]_{(\bm{0},\bm{0})}$.

  On the other hand as $I$ cuts out the affine diagonal, the set of homogeneous polynomials vanishing on $V(I)$ is exactly $I_\Delta$ so
  \[I[x_{\widehat{\tau}}^{-1}] \cap S[x_{\widehat{\tau}}^{-1}]_{(\bm{0},\bm{0})} = (I_{\Delta}[x_{\widehat{\tau}}^{-1}])_{(\bm{0},\bm{0})}.\]
  As $M[x_{\widehat{\tau}}^{-1}]$ is generated in degree $(\bm{0},\bm{0})$ by multiplying the above by $M[x_{\widehat{\tau}}^{-1}]$ we get the following:
  \[(IM[x_{\widehat{\tau}}^{-1}])_0 \cap M[x_{\widehat{\tau}}^{-1}]_{(0,0)} = (I_{\Delta}M[x_{\widehat{\tau}}^{-1}])_{(\bm{0},\bm{0})} = 0.\]
  This implies
  \[\Gamma(U_\tau,\varphi_*(\cO_X)) \cong  M[x_{\widehat{\tau}}^{-1}]_{(\bm{0},\bm{0})} = \Gamma(U_\tau,M^{\sim}).\]

  Moreover, the gluing maps are also isomorphic as the above identifications over the affine charts is induced from the $\Cl(X)$-graded map of $S$-modules $M\rightarrow R \cong S/I$.
\end{proof}

As a consequence of this lemma, the sheaf associated to a module with the appropriate support generated by monomials antisymmetric in $x$ and $y$ is the structure sheaf of the diagonal. This gives the following corollary which will be used to prove the sheaf statement in Theorem~\ref{thm:linear}.
\begin{corollary}
  \label{cor:lawrence-diagonal}
  If $X$ is a simplicial toric variety and if $M\subseteq \field[x_1^\pm,\ldots,x_n^\pm,y_1^\pm,\ldots,y_n^\pm]/I_{\Lambda(L)}$ is a $S$-submodule generated by monomials of the form $\frac{x^{\bm{u}}}{y^{\bm{u}}}$ for $\bm{u}\in \ZZ^n$ such that $1\in M$, then
  \[M^\sim \cong \cO_{\Delta}.\]
\end{corollary}
\begin{proof}
  Let $M^{\Pic}$ be the submodule of $M$ which is generated by all Picard degrees of $M$.
  That is to say the submodule of $M$ generated by all $m\in M_{\alpha}$ where $\alpha\in \Pic(X\times X)$.
  Set $M' = M/M^{\Pic}$ and Set $B$ to be the irrelevant ideal of $X\times X$. By \cite{CLS11}*{Exercise~5.3.5},
  to check that $(M')^{\sim} = 0$, it suffices to check that $(B^{\ell}M')_{\alpha} = 0$ for all $\ell\gg 0$
  and $\alpha\in \Pic(X\times X)$. By construction, $M'$ is zero in degrees in $\Pic(X\times X)$, and so $(M')^{\sim}=0$.

  To apply Lemma~\ref{lem:homogenize-diagonal}, let $I=\left<x_i-y_i \,|\, i=1,\ldots,n\right>$.
  As $1\in M^{\Pic}$ and $\operatorname{ann} M^{\Pic} = I_{\Lambda(L)}\subseteq I$,
  there is an injective map $S/I \rightarrow M^{\Pic}/IM^{\Pic}$.
  Moreover $M^{\Pic}$ is generated by Laurent monomials of the form  $\frac{x^{\bm{u}}}{y^{\bm{v}}}$
  for $\bm{u},\bm{v}\in \ZZ^n$ satisfying $\bm{v}\leq \bm{u}$. In the quotient $M^{\Pic}/IM^{\Pic}$,
  all of the Laurent monomials of that form are identified with the image of monomials of $S/I$
  and thus the map is surjective and thus is an isomorphism.
  Lemma~\ref{lem:homogenize-diagonal} then implies that $(M^{\Pic})^{\sim} \cong \cO_\Delta$ and thus $M^{\sim} \cong \cO_{\Delta}$.
  
\end{proof}

Another approach to the above corollary would be to work over the Picard graded subring of $S$ as in \cite{cox-paper}*{Theorem~3.11}, however, as everything else in this paper is in terms of the $\Cl(X)$-graded ring $S$, some translation would be required.

\subsection*{Linear Resolutions of the Diagonal}

The resolutions that the subdivision procedure constructs are essentially never minimal and so cannot be linear.
However, the existence of a linear minimal free resolution can often be detected from a nonminimal resolution.
In particular, if a module $M$ is resolved by a nonminimal resolution $F$ such that $F_0$ is generated in a single standard graded degree,
and all matrices representing the differentials in the resolution have only constant and linear entries,
then it follows that the minimal free resolution of $M$ is linear.

The modules arising from these barycentric subdivisions are always generated in degree $0$ in the standard grading.
However, even in the case of the diagonal of $\mathbb{P}^n$, when $\psi$ is the ceiling stratification, the complex $F_{\psi^{BC}}$ contains quadratic terms.

Nonetheless, with careful choice of stratification $\psi$ followed by a coarsening of the cell structure, we can get a complex with only linear and constant terms.
From this complex we can infer the existence of a linear resolution. The condition on the stratification which achieves this will be most easily phrased in terms of the quiver associated to the stratification.

\begin{definition}
  A quiver $Q$ labeled by a polynomial ring is \emph{linear} if all of the labels are variables.
\end{definition}

Then from this condition, a criterion for existence of linear resolutions becomes easy. Observe also that $Q_\psi$ is linear if and only if $\widetilde{Q}_\psi$ is linear.

\begin{proposition}
  \label{prop:linear}
  Let $\psi$ be a shifted ceiling stratification.
  If the quiver $\widetilde{Q}_\psi$ is linear, then the module $\pi_{\Lambda(L)}(M_{\psi^{BC}})$ has a linear resolution over $S$.
  Moreover, if the quiver $\widetilde{Q}_\psi$ has an edge $\bm{a}\rightarrow \bm{b}$ for which there does not exist a path from $\bm{a}$ to $\bm{b}$ consisting only of linear edges, then $\pi_{\Lambda(L)}(M_{\psi^{BC}})$ does not have a linear resolution.

  In particular as $\Mepsilon=\psi_{\Lambda(L)}(M_{\psi_{\epsilon}^{BC}})$, the above applies to $\Mepsilon$ with the quiver $\widetilde{Q}_{\psi_\epsilon}$.
\end{proposition}
\begin{proof}
  First, for the case of $\Mepsilon$, we use Proposition~\ref{prop:resolution} to reduce to the general case.
  
  We begin by showing that $\pi_{\Lambda(L)}(M_{\psi^{BC}})$ has a linear resolution over $S$. Let us define a cell structure on $L_{\RR}$ where the cells are in bijection with pairs of cells $(\sigma,\sigma')$ in the original cell complex on $L_{\RR}$ satisfying $\sigma\preceq \sigma'$. Define the closed cell $U_{(\sigma,\sigma')}$ to be the convex hull given by
  \[U_{(\sigma,\sigma')} = \hull\{v_\tau \,|\, \sigma\preceq \tau \preceq \sigma'\}.\]

  First, observe that each cell $\hull\{v_{\sigma_0},\ldots,v_{\sigma_d}\}$ of $L_{\RR}^{BC}$ is contained in $U_{(\sigma,\sigma')}$ if and only if $\sigma\preceq \sigma_0$ and $\sigma_d\preceq \sigma'$. Second, the dimension of $U_{(\sigma,\sigma')}$ is $\dim \sigma' - \dim \sigma$. The facets of $U_{(\sigma,\sigma')}$ are the cells $U_{(\sigma,\tau')}$ and $U_{(\tau,\sigma')}$ where $\sigma$ is a facet of $\tau$ and $\tau'$ is a facet of $\sigma'$.

  Let the cell structure given by this be denoted $L_{\RR}^{Quiver}$. As before, this induces a cell structure on $\Lambda(L)_{\RR}$ which we denote $\Lambda(L)_{\RR}^{Quiver}$. One can verify that the stratification function $\psi^{BC}$ remains constant on the open cells of $\Lambda(L)_{\RR}^{Quiver}$. Denote the resulting stratification by $\psi^{Quiver}$. This cell structure is finer than the original polyhedral cell structure given by $\psi$ and coarser than the barycentric subdivision. As $\Lambda(L)_{\RR}^{Quiver}$ is a coarsening of $\Lambda(L)_{\RR}^{BC}$, we get an injective map of complexes $F_{\psi^{Quiver}}\rightarrow F_{\psi^{BC}}$ with the map on chain complexes given by
  \[U_{(\sigma,\sigma')} \mapsto \sum_{\tau\subseteq L_{\RR}^{BC}, \tau\subseteq U_{(\sigma,\sigma')}} \tau.\]
  As both the underlying topological space and the stratification function are unchanged, this map of complexes induces an isomorphism on homology.
  
  As the module is generated in a single degree, it suffices to show that the complex $F_{\psi^{Quiver}}$ has differentials containing only constant and linear entries. For this purpose we carefully describe the differentials of this complex up to sign.
  \[\partial(U_{(\sigma,\sigma')}) = \sum_{\sigma\prec \tau,\dim \sigma + 1 = \dim \tau} \pm y^{\psi(\tau)-\psi(\sigma)}U_{(\tau,\sigma')} + \sum_{\tau'\prec \sigma', \dim \tau' + 1 = \dim \sigma'} \pm x^{\psi(\sigma')-\psi(\tau')}U_{(\sigma,\tau')}\]
  For the first summation either $\psi(\sigma)=\psi(\tau)$ or the edge $\psi(\tau)\rightarrow \psi(\sigma)$ is an edge in the quiver. As the quiver is linear, $y^{\psi(\tau)-\psi(\sigma)}$ is linear or constant. Similarly in the second summation we find that $x^{\psi(\sigma')-\psi(\tau')}$ is linear or constant. Thus the complex $F_{\psi^{Quiver}}$ contains only linear and constant terms, and consequently, the same holds for $\pi_{\Lambda(L)}(F_{\psi^{Quiver}})$.

  Finally, as the module $\pi_{\Lambda(L)}(F_{\psi^{Quiver},0})=\pi_{\Lambda(L)}(F_{\psi^{BC},0})$ is generated in a single standard graded degree and $H_0(F_{\psi^{Quiver}})=H_0(F_{\psi^{BC}})=M_{\psi^{BC}}$, the module $\pi_{\Lambda(L)}(M_{\psi^{BC}})$ has a linear resolution.

  For the partial converse, suppose that $\sigma\prec \sigma'$ is a pair of cells such that there is an edge $\psi(\sigma') \rightarrow \psi(\sigma)$ in the quiver $Q_\psi$ with a monomial label of degree $d\geq 2$ such that there does not exist a (directed) path consisting of linear edges from $\psi(\sigma')$ to $\psi(\sigma)$. Without loss of generality, we may assume that all paths from $\psi(\sigma')$ to $\psi(\sigma)$ are length $1$.
  We will show that $\beta_{1,d}(\pi_{\Lambda(L)}(M_{\psi^{BC}}))$ is nonzero and since $\pi_{\Lambda(L)}(M_{\psi^{BC}})$ is generated in degree $0$, this implies that $\pi_{\Lambda(L)}(M_{\psi^{BC}})$ does not have a linear resolution.

  Generalizing \cite{BS98}*{Theorem~1.11} to compatible $\ZZ^n$-stratifications, to show that $\beta_{1,d}$ is nonzero it suffices to find a $\bm{u}$ where
  \[\dim H_1\left((\Lambda(L)_{\RR}^{Quiver})_{\leq \bm{u}},(\Lambda(L)_{\RR}^{Quiver})_{<\bm{u}}\right)\neq 0.\]  
  See also \cite{FS24}*{Theorem~4.13 and Remark~4.14} as well as \cite{stratification-paper}*{Remark~2.16}.
  
  Fix $\bm{u}=(\psi(\sigma'),-\psi(\sigma))\in \ZZ^{2n}$
  and choose any path from $v_{\sigma'}$ to $v_{\sigma}$ through the 1-skeleton of $(\Lambda(L)_{\RR}^{Quiver})_{\leq \bm{u}}$.
  Such a path exists in $(\Lambda(L)_{\RR}^{Quiver})_{\leq \bm{u}}$ since it is contractible.
  This path gives a cycle in the relative homology $H_1\left((\Lambda(L)_{\RR}^{Quiver})_{\leq \bm{u}},(\Lambda(L)_{\RR}^{Quiver})_{< \bm{u}}\right)$
  as both $v_{\sigma'}$ and $v_{\sigma}$ are in $(\Lambda(L)_{\RR}^{Quiver})_{< \bm{u}}$.
  To show that $\beta_{1,d}(\pi_{\Lambda(L)}(M_{\psi^{BC}}))$ is nonzero it suffices to show that this cycle is nonzero.
  It suffices to show that no path through the 1-skeleton from $v_{\sigma'}$ to $v_{\sigma}$ is contained in $(\Lambda(L)_{\RR}^{Quiver})_{<\bm{u}}$.

  Let $v_{\tau}$ be any vertex of $(\Lambda(L)_{\RR}^{Quiver})_{\leq \bm{u}}$. First observe that $\psi(\sigma')\geq \psi(\tau)\geq \psi(\sigma)$ and so Proposition~\ref{prop:degree-quiver} implies there exists a path
  \[\psi(\sigma')\rightarrow \psi(\tau)\rightarrow \psi(\sigma).\]
  However, as every path from $\psi(\sigma')$ to $\psi(\sigma)$ is length 1, this implies $\psi(\sigma')=\psi(\tau)$ or $\psi(\tau)=\psi(\sigma)$. Thus every vertex in $(\Lambda(L)_{\RR}^{Quiver})_{\leq \bm{u}}$ has weight either $\psi^{BC}(v_{\sigma'})$ or $\psi^{BC}(v_\sigma)$. But this implies that any path through the 1-skeleton from $v_{\sigma'}$ to $v_{\sigma}$ must have an edge whose vertices have weight $\psi^{BC}(v_{\sigma'})$ and $\psi^{BC}(v_\sigma)$. But $\bm{u} = \mathbf{max}(\psi^{BC}(v_{\sigma'}),\psi^{BC}(v_\sigma))$ and so this edge is not in $(\Lambda(L)_{\RR}^{Quiver})_{< \bm{u}}$.

  Thus any path from $v_{\sigma'}$ to $v_{\sigma}$ gives a nonvanishing homology class in
  \[H_1\left((\Lambda(L)_{\RR}^{Quiver})_{\leq \bm{u}},(\Lambda(L)_{\RR}^{Quiver})_{<\bm{u}}\right).\]
  This gives $\beta_{1,d}(\pi_{\Lambda(L)}(M_{\psi^{BC}}))\neq 0$ as desired.
\end{proof}

Proposition~\ref{prop:linear} reduces our problem to giving stratifications where the quiver $Q_\psi$ is linear.
Unfortunately, using only the ceiling stratification, the quiver $Q_\psi$ need not be linear and consequently, $\pi_{\Lambda(L)}(M_{\psi^{BC}})$ need not have a linear resolution. While many small examples may in fact give linear resolutions by chance, the following example exhibits a nonlinear minimal free resolution and demonstrates that the $\epsilon$-shift is necessary to get linear resolutions in general.
\begin{example}
  \label{ex:nonlinear}

  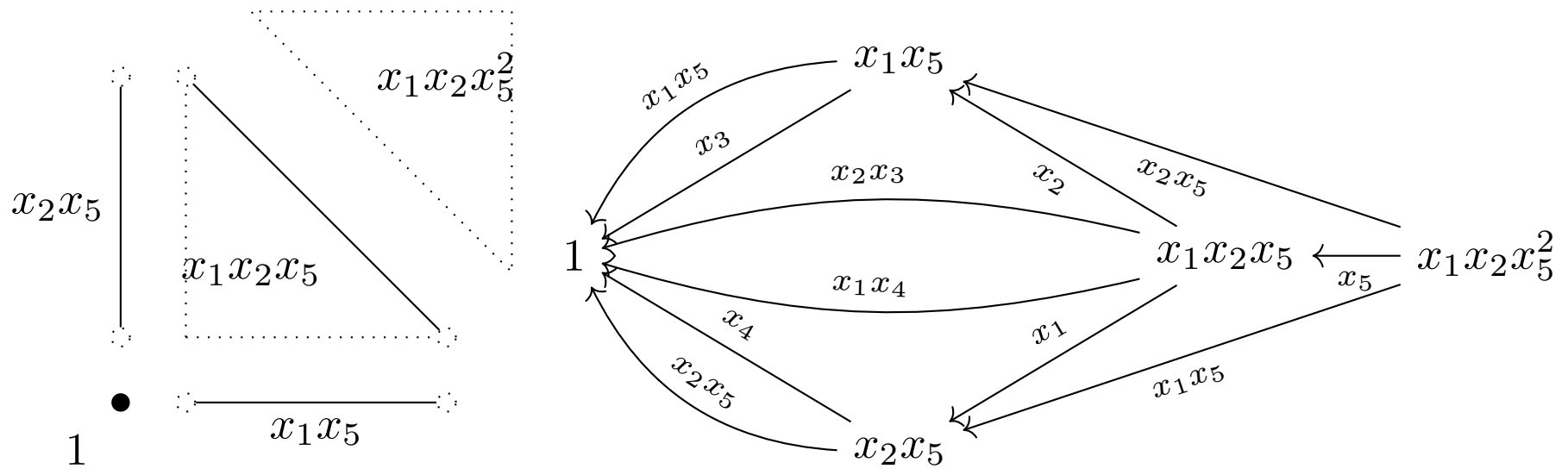
\begin{figure}[b]
    \begin{tikzpicture}[scale=2]

      \node[circle, fill, inner sep=0.05cm] (or) at (0,0) {};
      \node[below left = 0.1cm of or] {1};
      
      \begin{scope}[xshift=0.25cm]
        \node[circle,draw, dotted,inner sep=0.05cm] (a) at (0,0) {};
        \node[circle,draw, dotted,inner sep=0.05cm] (b) at (1,0) {};
        \draw (a) -- node[below] {$x_1x_5$} (b);
      \end{scope}

      \begin{scope}[yshift=0.25cm]
        \node[circle,draw, dotted,inner sep=0.05cm] (a) at (0,0) {};
        \node[circle,draw, dotted,inner sep=0.05cm] (b) at (0,1) {};
        \draw (a) -- node[left] {$x_2x_5$} (b);
      \end{scope}

      \begin{scope}[xshift=0.25cm,yshift=0.25cm]
        \coordinate (a) at (0,0);
        \node[circle,draw, dotted,inner sep=0.05cm] (b) at (1,0) {};
        \node[circle,draw, dotted,inner sep=0.05cm] (c) at (0,1) {};
        \draw[dotted] (a) -- (b);
        \draw[dotted] (a) -- (c);
        \draw (b) -- (c);
        \node at (0.25,0.25) {$x_1x_2x_5$};
      \end{scope}
      
      \begin{scope}[xshift=0.5cm,yshift=0.5cm]
        \coordinate (a) at (1,1);
        \coordinate (b) at (1,0);
        \coordinate (c) at (0,1);
        \draw[dotted] (a) -- (b) -- (c) -- (a);
        \node at (0.75,0.75) {$x_1x_2x_5^2$};
      \end{scope}
      
    \end{tikzpicture}
    \begin{tikzpicture}
      \node (p) at (0,0) {$1$};
      \node (h) at (2.5,1.5) {$x_1x_5$};
      \node (d) at (5,0) {$x_1x_2x_5$};
      \node (v) at (2.5,-1.5) {$x_2x_5$};
      \node (r) at (7,0) {$x_1x_2x_5^2$};

      \draw[->] (h) -- node[midway, above, sloped] {$\scriptstyle x_3$} (p);
      \draw[->] (h) edge[bend right] node[midway, above, sloped] {$\scriptstyle x_1x_5$} (p);
      \draw[->] (v) -- node[midway, above, sloped] {$\scriptstyle x_4$} (p);
      \draw[->] (v) edge[bend left] node[midway, above, sloped] {$\scriptstyle x_2x_5$} (p);
      \draw[->] (d) edge[bend right=15] node[midway, above, sloped] {$\scriptstyle x_2x_3$} (p);
      \draw[->] (d) edge[bend left=15] node[midway, above, sloped] {$\scriptstyle x_1x_4$} (p);

      \draw[->] (d) -- node[midway, below, sloped] {$\scriptstyle x_2$} (h);
      \draw[->] (d) -- node[midway, above, sloped] {$\scriptstyle x_1$} (v);
      
      \draw[->] (r) -- node[midway, below, sloped] {$\scriptstyle x_1x_5$} (v);
      \draw[->] (r) -- node[midway, below, sloped] {$\scriptstyle x_2x_5$} (h);
      \draw[->] (r) -- node[midway, below, sloped] {$\scriptstyle x_5$} (d);

    \end{tikzpicture}
    \caption{The strata of the stratification ceiling stratification for $X$ from Example~\ref{ex:nonlinear} and the corresponding quiver.}
    \label{fig:strata-and-quiver}
  \end{figure}

  Let $X$ be the blowup of $\mathbb{P}^1\times \mathbb{P}^1$ at one torus-invariant point.
  Then the Cox ring for $X$ is $R=\field[x_1,x_2,x_3,x_4,x_5]$ where $x_5$ is the variable arising from the blowup.
  
  A basis for the lattice $M_X\subseteq \ZZ^5$ is given by $(1,0,-1,0,1)$ and $(0,1,0,-1,1)$ which we use as the $x$ and $y$ coordinate directions respectively in the depiction of the strata in Figures~\ref{fig:strata-and-quiver} and \ref{fig:es-quiver}. In particular, this basis is chosen such that the rays of the fan for $X$ depicted below give the normal vectors for the hyperplanes in $\RR^2$ cutting out the stratification.
  
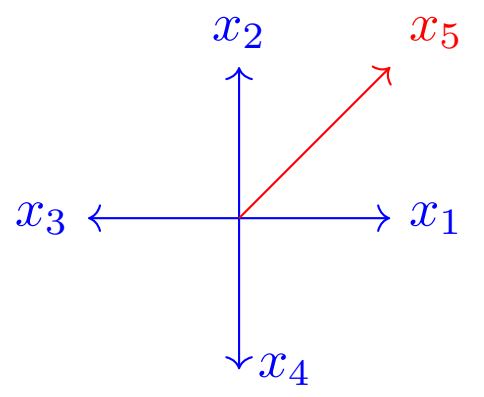
\begin{center}
    \begin{tikzpicture}
      \draw[blue,->] (0,0) -- (1,0) node[right] {$x_1$};
      \draw[blue,->] (0,0) -- (0,1) node[above] {$x_2$};
      \draw[blue,->] (0,0) -- (-1,0) node[left] {$x_3$};
      \draw[blue,->] (0,0) -- (0,-1) node[right] {$x_4$};
      \draw[red,->] (0,0) -- (1,1) node[above right] {$x_5$};
    \end{tikzpicture}
  \end{center}
  Using the ceiling stratification yields the stratification and quiver illustrated in Figure~\ref{fig:strata-and-quiver}.
  Then as there are quadratic edges which have no corresponding linear path, such the edge with monomial label $x_1x_5$ , Proposition~\ref{prop:linear} implies that the module $\pi_{\Lambda(L)}(M_{\psi^{BC}})$ arising from the ceiling stratification $\psi$ will not have a linear resolution. A word of caution, the diagram in Figure~\ref{fig:strata-and-quiver} is of $Q_\psi$ but the proposition is stated for $\widetilde{Q}_\psi$, so while there is a linear path in $Q_\psi$ from the vertex labeled $x_1x_5$ to the vertex labeled $1$, these are not between the same vertices in $\widetilde{Q}_\psi$.
  Taking the minimal free resolution yields the following Betti table in the standard $\ZZ$-grading.
  \[\begin{array}{c|ccc}
     & 0 & 1 & 2\\ \hline
    0 & 5 & 7 & 1  \\
    1 & - & 2 & 3 \\
  \end{array}\]

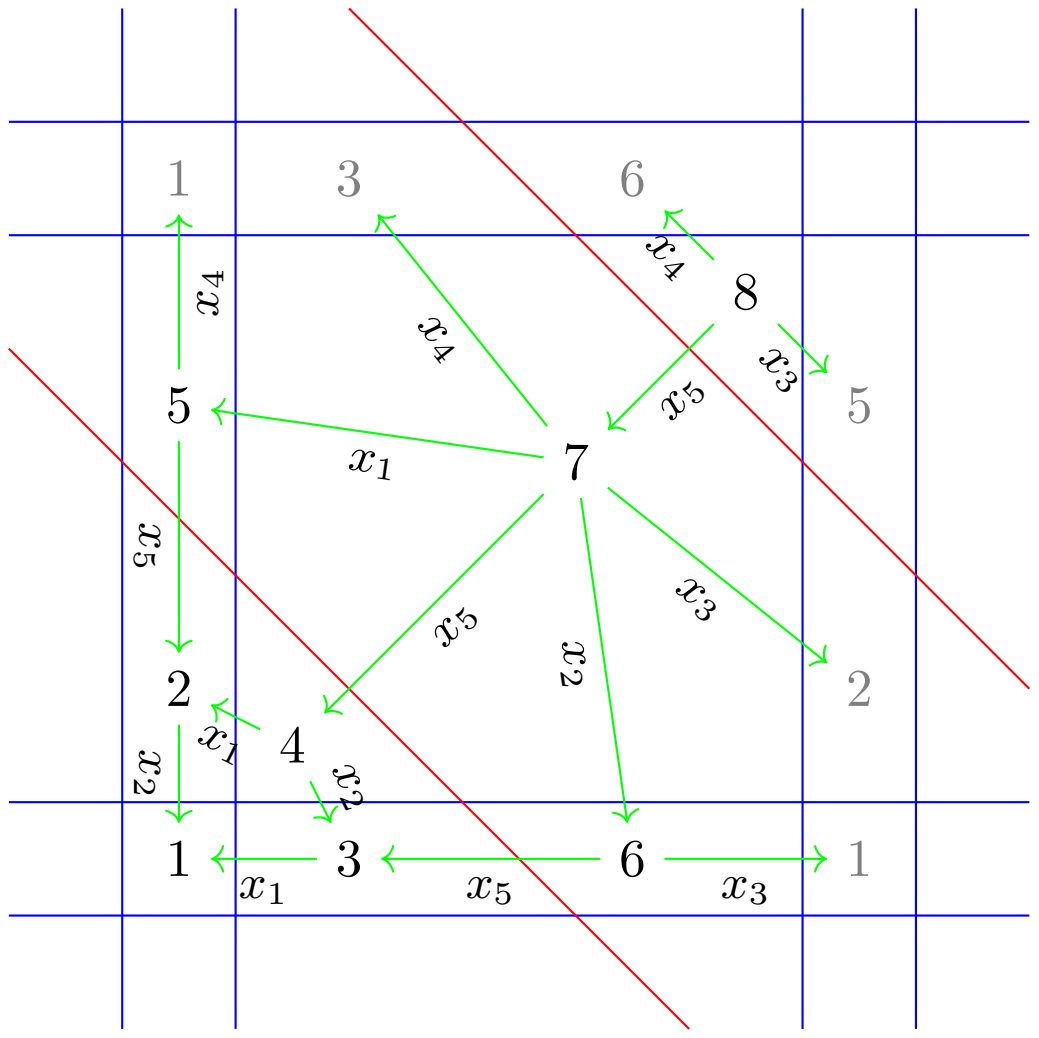
  \begin{figure}[b]
    \begin{center}
      \begin{tikzpicture}[scale=0.75]
        \coordinate  (origin) at (0,0);
        \coordinate (br) at (6,0);
        \coordinate (tr) at (6,6);
        \coordinate (tl) at (0,6);

        \coordinate (epsilon) at (1,1);

        \draw[blue] ($(origin)+(epsilon)-(0,2)$) -- ($(tl)+(epsilon)+(0,1)$);
        \draw[blue] ($(origin)+(epsilon)-(2,0)$) -- ($(br)+(epsilon)+(1,0)$);
        \draw[blue] ($(br)+(epsilon)-(0,2)$) -- ($(tr)+(epsilon)+(0,1)$);
        \draw[blue] ($(tl)+(epsilon)-(2,0)$) -- ($(tr)+(epsilon)+(1,0)$);

        \draw[blue] ($(origin)-(0,1)$) -- ($(tl)+(0,2)$);
        \draw[blue] ($(origin)-(1,0)$) -- ($(br)+(2,0)$);
        \draw[blue] ($(br)-(0,1)$) -- ($(tr)+(0,2)$);
        \draw[blue] ($(tl)-(1,0)$) -- ($(tr)+(2,0)$);

        \draw[red] ($(tl)+(epsilon)+(epsilon)$) -- ($(br)+(epsilon)+(epsilon)$);
        \draw[red] ($(origin)+(epsilon)+(epsilon)-(3,-3)$) -- ($(origin)+(epsilon)+(epsilon)+(3,-3)$);

        \node (v1) at (0.5,0.5) {1};
        \node (v2) at (0.5,2) {2};
        \node (v3) at (2,0.5) {3};
        \node (v4) at (1.5,1.5) {4};
        \node (v5) at (0.5,4.5) {5};
        \node (v6) at (4.5,0.5) {6};
        \node (v7) at (4,4) {7};
        \node (v8) at (5.5,5.5) {8};

        \begin{scope}[green]
          \node[gray] (v1') at ($(v1)+(tl)$) {$1$};
          \node[gray] (v1'') at ($(v1)+(br)$) {$1$};
          
          \node[gray] (v2') at ($(v2)+(br)$) {$2$};
          \node[gray] (v3') at ($(v3)+(tl)$) {$3$};
          \node[gray] (v5') at ($(v5)+(br)$) {$5$};
          \node[gray] (v6') at ($(v6)+(tl)$) {$6$};

          \draw[->] (v2) -- node[midway, below, sloped, black] {\small $x_2$} (v1);
          \draw[->] (v3) -- node[midway, below, sloped, black] {\small $x_1$} (v1);
          \draw[->] (v4) -- node[midway, below, sloped, black] {\small $x_1$} (v2);
          \draw[->] (v4) -- node[midway, above, sloped, black] {\small $x_2$} (v3);
          \draw[->] (v5) -- node[midway, below, sloped, black] {\small $x_5$} (v2);
          \draw[->] (v5) -- node[midway, below, sloped, black] {\small $x_4$} (v1');
          \draw[->] (v6) -- node[midway, below, sloped, black] {\small $x_5$} (v3);
          \draw[->] (v6) -- node[midway, below, sloped, black] {\small $x_3$} (v1'');
          \draw[->] (v7) -- node[midway, below, sloped, black] {\small $x_5$} (v4);
          \draw[->] (v7) -- node[midway, below, sloped, black] {\small $x_1$} (v5);
          \draw[->] (v7) -- node[midway, below, sloped, black] {\small $x_2$} (v6);
          \draw[->] (v7) -- node[midway, below, sloped, black] {\small $x_3$} (v2');
          \draw[->] (v7) -- node[midway, below, sloped, black] {\small $x_4$} (v3');
          \draw[->] (v8) -- node[midway, below, sloped, black] {\small $x_5$} (v7);
          \draw[->] (v8) -- node[midway, below, sloped, black] {\small $x_3$} (v5');
          \draw[->] (v8) -- node[midway, below, sloped, black] {\small $x_4$} (v6');
        \end{scope}

      \end{tikzpicture}
      \[\begin{array}{|l|c|c|c|c|c|c|c|c|}
      \hline
      \text{strata} & 1 & 2 & 3 & 4 & 5 & 6 & 7 & 8 \\\hline
      \text{monomial} & 1 & x_2 & x_1 & x_1x_2& x_2x_5 & x_1x_5& x_1x_2x_5 & x_1x_2x_5^2 \\\hline
      \end{array}\]
    \end{center}
    \caption{The stratification and quiver resulting from shifting the hyperplanes by $\epsilon=(1/3,1/3,0,0,2/3)$ with the edges of the quiver labeled by the monomials corresponding to those edges.}
    \label{fig:es-quiver}
  \end{figure}

  Theorem~\ref{thm:linear} suggests that if we take a generic $\epsilon$, the resulting module should have a linear resolution. Here we take $\epsilon=(1/3,1/3,0,0,2/3)$ and get the cell structure and quiver illustrated in Figure~\ref{fig:es-quiver}.
  In particular, the resulting quiver is now linear. In this case,
  \[\pi_{\Lambda(L)}(M_{\psi_\epsilon^{BC}}) = S\cdot \left\{1,\frac{x_2}{y_2},\frac{x_1}{y_1},\frac{x_1x_2}{y_1y_2},\frac{x_2x_5}{y_2y_5},\frac{x_1x_5}{y_1y_5},\frac{x_1x_2x_5}{y_1y_2y_5},\frac{x_1x_2x_5^2}{y_1y_2y_5^2}\right\} \subseteq \field[\ZZ^{2n}/\Lambda(L)].\]
  A quick Macaulay2 \cite{M2} computation shows that this module has a minimal Betti table given by
  \[\begin{array}{c|ccc}
    & 0 & 1 & 2 \\\hline
    0 & 8 & 16 & 8
  \end{array}\]
\end{example}

A careful observation of Example~\ref{ex:nonlinear} shows that the nonlinearity arises in large part as a consequence of the non-transversality of the intersection of the hyperplanes at the lattice points of $L$.
We say that a collection of hyperplanes in dimension $d$ has \emph{transverse intersections} if every intersection of $c$ hyperplanes for $c\leq d$ is codimension exactly $c$. Assuming that the stratification has a cell complex arising from a hyperplane arrangement with transverse intersections, the linearity does in fact follow. 

Moreover, as the hyperplane arrangement associated to a ceiling stratification contains finitely many equivalence classes of hyperplanes, to give transverse intersections we can simply offset each hyperplane by a generic shift. As we previously constructed in Section~\ref{sec:background}, we can shift our stratification by some $\epsilon\in \RR^{n}$ to get a new stratification by hyperplanes which we denote $\psi_\epsilon$.
To recall, for a point  $\bm{v}\in L_{\RR}$, the value of $\psi_{\epsilon}$ is given by
\[\psi_{\epsilon}(\bm{v}) = \lceil \bm{v}-\epsilon \rceil.\]
This construction will form the basis of the proof of the main theorem.

\begin{proof}[Proof of Theorem~\ref{thm:linear}]
  For $\Delta\subseteq X\times X$ the diagonal subvariety,
  the ideal for the diagonal is the toric ideal corresponding to the lattice
  $\Lambda(L)\subseteq \ZZ^{2n}$. Equivalently, $I_\Delta = I_{\Lambda(L)}$.
  Let $\psi = \psi_\epsilon$ be the shifted ceiling stratification corresponding to $\epsilon$.
  Then as the image of $\psi$ in $\ZZ^n$ is exactly all vectors of the form
  $\lceil\bm{p}-\epsilon\rceil$ for $\bm{p}\in L$, by Proposition~\ref{prop:resolution},
  the modules $\Mepsilon$ and $\pi_{\Lambda(L)}(M_{\psi^{BC}})$ are equal.
  Moreover, the same proposition along with the fact that $\dim (L_{\RR}) = \dim X$
  implies that $\pdim \Mepsilon \leq \dim X$.

  As $\epsilon_i<1$, it follows that $\psi_{\epsilon}(\bm{0}) = 0$ and thus
  Corollary~\ref{cor:lawrence-diagonal} implies that if $X$ is simplicial $(\Mepsilon)^{\sim} = \cO_\Delta$.
  Finally, as $\codim \Delta=\dim X$, we find that $\pdim \Mepsilon \geq \dim X$
  and thus $\pdim \Mepsilon = \dim X$.

  Next, we show the existence of a linear resolution.
  By Proposition~\ref{prop:linear}, it suffices to find a compatible $\ZZ^n$-stratification
  of $L_{\RR}$ that has a linear quiver. As the number of hyperplane classes defining
  the cell structure of $\psi$ is finite up to the $L$-action, for a generic $\epsilon$,
  the hyperplanes that define the chambers of the cell structure
  for the stratification $\psi$ intersect transversely. For the remainder of the proof,
  fix a choice of $\epsilon$ such that these intersections are transverse.

  Now let $\bm{a}\rightarrow \bm{b}$ be an edge in the quiver $\widetilde{Q}_{\psi}$. Then there exists $\sigma\prec \sigma'$ such that $\psi(\sigma)=\bm{b}$ and $\psi(\sigma')=\bm{a}$ and where $\sigma$ is a facet of $\sigma'$. Then as $\sigma$ is a facet of $\sigma'$ and the hyperplanes for the stratification are transverse, $\sigma$ is the intersection of $\sigma'$ with exactly one hyperplane. It follows that $|\bm{a}-\bm{b}|\leq 1$. Since $X$ is a complete normal toric variety, $M_X\cap \NN^{n}=\{0\}$. Thus by Lemma~\ref{lem:quiver-basics}, the quiver contains no cycles so $\bm{a}\neq \bm{b}$ and thus $|\bm{a}-\bm{b}|=1$. Thus the edge $\bm{a}\rightarrow \bm{b}$ is labeled by a variable and so the quiver is linear. Thus by Proposition~\ref{prop:linear}, the module $\pi_{\Lambda(L)}(M_{\psi^{BC}})$ has a linear resolution.
\end{proof}

Now we conclude with the proof of Corollary~\ref{cor:truncation}.
\begin{proof}[Proof of Corollary~\ref{cor:truncation}]
  Fix $\epsilon\in \RR^n$ such that $\Mepsilon$ has a linear resolution of length $\dim X$
  and fix an ample degree $\alpha\in \Cl(X)$ with $\alpha\gg \bm{0}$.
  Recall that $S$ is given a $\Cl(X)$-grading by coarsening the $\Cl(X\times X)$ grading.

  The module $\Mepsilon$ is a finitely generated $\Cl(X)$-graded submodule of $\field[\ZZ^{2n}/\Lambda(L)]$.
  In particular, there exists a finite collection of monomials in $\field[\ZZ^{2n}]$ such that
  their image in $\field[\ZZ^{2n}/\Lambda(L)]$ generates $\Mepsilon$.
  By clearing denominators, there exists a monomial $m\in S$ $m\cdot \Mepsilon$
  is a submodule of $\field[\NN^{2n}/\Lambda(L)]\cong S/I_{\Delta}$.

  Then let $\alpha$ be any degree such that $\alpha-\deg m$ is effective,
  there exists a monomial $m'\in S$ with $\deg m' = \alpha$ where
  $m'\cdot \Mepsilon$ is a submodule of $S/I_{\Delta}$. As $\Mepsilon$ is generated
  in degree $\bm{0}\in \Cl(X)$, the module $m'\cdot \Mepsilon$ is generated
  in degree $\alpha$ and has a linear resolution of length $\dim X$.
  Thus for any sufficiently ample degree, there exists a submodule of $S/I_{\Delta}$
  generated in that degree with linear resolution of length $\dim X$.
\end{proof}


\section{Dual Complexes for Stratifications}
\label{sec:dual}

While we are not in general able to give a minimal free resolution of $\Mepsilon$, by a partial minimization we can construct a cellular resolution of $\Mepsilon$ with a particularly simple structure. As described in Definition~\ref{def:dual-complex}, there exists a simplicial complex $\Dual_\psi$ supporting this partially minimized resolution.

The name ``dual complex'' for $\Dual_\psi$ is justified by observing that the vertices of the complex are in bijective correspondence with the strata of the stratification. Moreover, each edge in the quiver $\widetilde{Q}_\psi$ contributes one edge of the $1$-skeleton. However, not every edge in this dual complex is an edge of the quiver.

To construct a resolution of the diagonal supported on the dual complex, we will show that there exists a homogeneous acyclic matching on the cell structure of $\Lambda(L)_{\RR}^{BC}$ that minimizes the resolution to one supported on the dual complex of $\psi$.

\subsection*{Acyclic Matchings on the Barycentric Subdivision}

We begin by recalling some background on acyclic matchings. A good resource for this is Batzies' PhD thesis \cite{Batzies}. There is an article by the same name by Batzies and Welker, but the thesis is more comprehensive and so is the one we cite here.
For a regular CW complex, the cells form a graph with edges $\tau \rightarrow \sigma$ if $\sigma$ is a facet of $\tau$. For a CW complex $\mathcal{X}$, denote this graph by $\graph(\mathcal{X})$. A matching $B$ on a CW complex is a collection of edges from the graph such that every cell occurs at most once. A matching is \emph{acyclic} if the graph obtained by reversing all edges in $B$ is acyclic. Denote this reversed graph by $\graph^{B}(\mathcal{X})$. The \emph{critical cells} of an acyclic matching are the cells that do not appear in any edge of the acyclic matching. To be suggestive of this reversal, we will write all of the edges in the acyclic matching in reverse. That is, we will write $\sigma \leftarrow \tau$ instead of $\tau\rightarrow \sigma$ for edges in the acyclic matching. As our cell complexes are not finite, some care is required. We need the acyclic matchings to be proper as defined in \cite{Batzies}*{Definition~3.2.10}. As all of our CW complexes are locally finite, it suffices to show that the paths in $\graph^{B}(\mathcal{X})$ are all finite, as this will imply that the universal grading defined in \cite{Batzies}*{Definition~3.2.5} is compact.

Discrete Morse theory gives an explicit homotopy equivalence between the original CW complex and a new CW complex whose cells are in bijection with the critical cells of the acyclic matching \cite{Batzies}*{Theorem~3.2.9 and Definition~3.2.13}.
We will denote this homotopy equivalence by $\homotopy(B)$ and for a CW complex $\mathcal{X}$ we will denote the resulting complex $\homotopy(B)(\mathcal{X})$.
Moreover, for a pair of critical cells $\tau,\sigma$ satisfying $\dim \tau = \dim \sigma - 1$, the cell $\tau$ is in the boundary of $\sigma$ if there exists a path in $\graph^{B}(\mathcal{X})$ from $\sigma$ to $\tau$.

In the context of cellular resolutions, if the monomial labels associated to cells $\sigma$ and $\tau$ agree for all edges $\tau\leftarrow \sigma$ in the acyclic matching, then we say that the matching is \emph{homogeneous}. In this circumstance, the complex resulting from the acyclic matching also has a labeling by monomials which resolves the same module \cite{Batzies}*{Theorem~3.3.1}. This extends immediately to a description in terms of compatible $\ZZ^n$-stratifications.

With the background on acyclic matchings complete, we move to defining some terminology for chains of cells. These will be helpful in defining the acyclic matching that yields the dual complex.

\begin{definition}
  Given a chain of cells $\sigma_0\prec \cdots \prec \sigma_d$ associated to a compatible $\ZZ^n$-stratification, we say the chain is
\begin{enumerate}
\item \emph{degree-strict} if the value of $\psi$ is distinct on each cell, and
\item \emph{degree-maximal} if the chain is degree-strict, $\sigma_d$ is the maximal cell where $\psi$ takes the value $\psi(\sigma_d)$, and for $0\leq i<d$, the cell $\sigma_i$ is maximal among cells where $\psi$ takes the value $\psi(\sigma_i)$ such that $\sigma_i\prec \sigma_{i+1}$. 
\end{enumerate}

Finally, define the \emph{degree-maximal tail} to be the longest contiguous degree-maximal chain ending at $\sigma_d$.

\end{definition}

We begin with a quick lemma about degree-maximal subchains.

\begin{lemma}
  \label{lem:degree-maximal-subchains}
  If $\sigma_0\prec \cdots \prec \sigma_d$ is a degree-maximal chain associated to a shifted ceiling stratification, then every subchain containing $\sigma_d$ is degree-maximal.
\end{lemma}
\begin{proof}
  It suffices to show that $\sigma_0\prec \cdots \prec \widehat{\sigma}_s\prec \cdots \prec \sigma_d$ is degree-maximal for $s<d$.
  In particular, it suffices to show that $\sigma_{s-1}$ is maximal among cells $\tau$ satisfying $\tau\prec \sigma_{s+1}$ and $\psi(\tau)=\psi(\sigma_{s-1})$.
  So let $\tau$ be a cell satisfying $\tau \prec \sigma_{s+1}$ and $\psi(\tau)=\psi(\sigma_{s-1})$.

  Set $\bm{u} = \psi(\sigma_s)$. As both $\overline{\sigma_{s+1}}$ and $(L_{\RR})_{\leq \bm{u}}$ are polyhedral their intersection is polyhedral.
  Moreover, $(L_{\RR})_{\leq \bm{u}}$ is a closed union of cells from $L_{\RR}$ and so the intersection
  $\overline{\sigma_{s+1}}\cap (L_{\RR})_{\leq \bm{u}}$ is the closed union of faces of $\sigma_{s+1}$.
  However as the intersection is also polyhedral, this implies that the intersection is the closure of a single face of $\sigma_{s+1}$.
  As the original chain is degree-maximal, $\sigma_s$ is the unique maximal cell of this intersection.
  Moreover, $\psi(\tau)\leq \bm{u}$ and so $\tau\subseteq \overline{\sigma_{s+1}}\cap (L_{\RR})_{\leq \bm{u}}=\overline{\sigma_s}$.
  Thus $\tau$ is a face of $\sigma_s$.
  As the original chain is degree-maximal and $\psi(\tau)=\psi(\sigma_{s-1})$, this implies that $\tau\preceq \sigma_{s-1}$.
  Thus $\sigma_{s-1}$ is maximal among faces $\tau$ of $\sigma_{s+1}$ with $\psi(\tau)=\psi(\sigma_{s-1})$
  and so the chain $\sigma_0\prec \cdots \prec \widehat{\sigma}_s\prec \cdots \prec \sigma_d$ is degree-maximal.
\end{proof}

Using these properties of chains, we can define an acyclic matching on the cell structure of $L_{\RR}^{BC}$.

\begin{definition}
  \label{def:quiver-matching}
  Let $\psi$ be a compatible $\ZZ^n$-stratification on $L_{\RR}$.
  For a chain of cells $\sigma_0\prec \cdots \prec \sigma_d$ let $0<s\leq d$ be the unique integer, if it exists, such that
  \begin{enumerate}
  \item $\psi(\sigma_{s-1}) = \psi(\sigma_s)$, and
  \item $\sigma_s\prec \cdots \prec \sigma_d$ is a degree-maximal chain.
  \end{enumerate}
  Associate to the chain the following edge in the graph $\graph(L_{\RR}^{BC})$:
  \[\hull\{v_{\sigma_0},\ldots,\widehat{v}_{\sigma_s},\ldots,v_{\sigma_d}\}\leftarrow \hull\{v_{\sigma_0},\ldots,v_{\sigma_s},\ldots,v_{\sigma_{d}}\}.\]

  The \emph{quiver matching associated to} $\psi$, denoted $\Quiv_\psi$, is the set of edges of $\graph(L_{\RR}^{BC})$ associated to such chains of cells.
\end{definition}

As an abuse of notation, we will say that a cell $\hull \{v_{\sigma_0},\ldots,v_{\sigma_d}\}\subseteq \Lambda(L)_{\RR}^{BC}$ is \emph{degree-strict} or \emph{degree-maximal} if the corresponding chain $\sigma_0\prec \cdots \prec \sigma_d$ is degree-strict or degree-maximal respectively. Moreover, for a cell $\alpha = \hull \{v_{\sigma_0},\ldots,v_{\sigma_d}\}$, the \emph{degree sequence associated to $\alpha$} is the chain $\psi(\sigma_0)\leq \cdots \leq \psi(\sigma_d)$.

\begin{proposition}
  Let $L\subseteq \ZZ^n$ be a lattice with $L \cap \NN^n=\{0\}$. If $\psi$ is a shifted ceiling stratification on $L$ then the quiver matching $\Quiv_\psi$ is a homogeneous proper acyclic matching.
\end{proposition}
\begin{proof}
  That the matching is homogeneous is immediate.

  We begin by showing that for any $\alpha$ there is at most one cell $\beta$ such that $\alpha \leftarrow \beta$ is an edge in $\Quiv_\psi$.
  Index the vertices of $\alpha$ as
  \[\alpha = \hull\{v_{\sigma_0},\ldots,v_{\sigma_{s-1}},v_{\sigma_{s+1}},\ldots,v_{\sigma_d}\}\]
  such that $\sigma_0\prec \cdots \prec \sigma_{s-1}\prec \sigma_{s+1}\prec \cdots \prec \sigma_d$
  and such that the chain $\sigma_{s+1}\prec \cdots\prec \sigma_d$ is the degree-maximal tail.
  Observe first that if $\alpha\leftarrow \beta$ is an edge in $\Quiv_\psi$ then the chain associated to it cannot be degree-maximal, and so $s\neq 0$. Moreover, the degree maximal tail of the chain associated to $\beta$ must contain the degree maximal chain $\sigma_{s+1}\prec \cdots \prec \sigma_d$ and so it must be that $\beta = \hull\{v_{\sigma_0},\ldots,v_{\sigma_{s-1}},v_{\sigma_s},v_{\sigma_{s+1}},\ldots,v_{\sigma_d}\}$ for some cell $\sigma_s$ satisfying $\sigma_{s-1}\prec \sigma_s \prec \sigma_{s+1}$. Then as $\alpha\leftarrow \beta$ is an edge in $\Quiv_\psi$, it follows that $\sigma_s$ is maximal among faces of $\sigma_{s+1}$ satisfying $\psi(\sigma_{s-1})=\psi(\sigma_s)$. By Lemma~\ref{lem:unique-maximal-cell}, the choice of $\sigma_s$ is unique and so $\beta$ is unique. 

  Next, to show that it is a matching, we need to show that for $\beta=\hull\{v_{\sigma_{0}},\ldots,v_{\sigma_d}\}$, there exist at most one edge of the form $\alpha \leftarrow \beta$ in $\Quiv_\psi$. Such an edge, if it exists, amounts to an appropriate choice of $s$ with $0<s\leq d$. As the chain $\sigma_{s}\prec \cdots \prec \sigma_d$ is degree-maximal, $s$ must be the maximal integer satisfying $\psi(\sigma_{s-1})=\psi(\sigma_s)$, so the choice is unique.

  Finally, suppose that $\alpha_1\rightarrow \beta_1 \rightarrow \alpha_2 \rightarrow \cdots$ is a potentially infinite path in $\graph^{\Quiv_\psi}(\Lambda(L)_{\RR}^{BC})$ such that $\alpha_i\leftarrow \beta_i$ is an edge in $\Quiv_\psi$ for $i=1,2,\ldots$.
  To show that the matching is acyclic and proper, it suffices to show that the path is finite.
  The matching is homogeneous so the value of $\psi^{BC}$ is non-increasing in any such path.
  Then, as the poset $\img \psi^{BC}\subseteq \ZZ^{2n}$ is well-founded,
  it suffices to consider chains where $\psi^{BC}(\alpha_i)=\psi^{BC}(\alpha_{i+1})$ for all $i$.

  Fix an index $i$. Then there exist cells $\sigma_0,\ldots,\sigma_d\subseteq L_{\RR}$ and integers $0\leq r,s\leq d$ such that
  \begin{align*}
    \alpha_i&=\hull\{v_{\sigma_0},\ldots,v_{\sigma_{r-1}},v_{\sigma_{r+1}},\ldots,v_{\sigma_d}\},\\
    \beta_i &= \hull\{v_{\sigma_0},\ldots,v_{\sigma_{r}},\ldots,v_{\sigma_d}\},\\
    \alpha_{i+1} &= \hull\{v_{\sigma_0},\ldots,v_{\sigma_{s-1}},v_{\sigma_{s+1}},\ldots,v_{\sigma_d}\}.\\
  \end{align*}
  
  For $\beta_{i+1}$ there are a few possibilities to consider.
  First $s\neq d$ as $\psi^{BC}(\alpha_i)=\psi^{BC}(\alpha_{i+1})$.
  Moreover if $s\neq r-1$ then the chain $\sigma_{r}\prec \cdots \prec \sigma_d$ (omitting $\sigma_s$ if $s>r$)
  is degree-maximal either immediately or by Lemma~\ref{lem:degree-maximal-subchains}.
  Moreover, it is the degree-maximal tail of $\sigma_0\prec \cdots \prec \widehat{\sigma}_s\prec \cdots \prec \sigma_d$. Thus if $\alpha_{i+1}\leftarrow \beta_{i+1}$ is an edge in $\Quiv_\psi$, then there exists a cell $\tau\subseteq L_{\RR}$ such that
  \[\beta_{i+1} = \hull\{v_{\sigma_0},\ldots,v_{\sigma_{r-1}},v_{\tau},v_{\sigma_{r}},\ldots,v_{\sigma_d}\}.\]
  Again in the above equation, we omit $\sigma_s$ from the appropriate location.
  However, $\psi(\sigma_{r-1})=\psi(\sigma_r)$ and thus $\psi(\tau)=\psi(\sigma_r)$ and so the chain $\tau\prec \sigma_r\prec \cdots \prec \sigma_d$ cannot be degree-maximal. Thus $\alpha_{i+1}\leftarrow \beta_{i+1}$ cannot be an edge in $\Quiv_\psi$.

  If $s = r-1$, observe that as the chain $\sigma_{r}\prec \cdots \prec \sigma_d$ is degree-maximal and the corresponding vertices appear in $\alpha_{i+1}$, the degree-maximal tail arising from the vertices of $\beta_{i+1}$ must contain this chain as a subchain and thus must be strictly longer. In particular the length of the degree-maximal tail for the chain corresponding to $\beta_{i+1}$ is longer than the degree-maximal tail for the chain corresponding to $\beta_i$. As the length of the full chain does not change, the path $\alpha_1\rightarrow \beta_1\rightarrow \alpha_2 \rightarrow \beta_2 \rightarrow \cdots$ in $\graph^{\Quiv_\psi}(L_{\RR}^{BC})$ terminates.
  Thus the matching is acyclic.
\end{proof}

To describe the resulting CW complex after applying the homotopy associated to the matching, we must first identify the critical cells.

\begin{lemma}
  \label{lem:critical-cells}
  The critical cells of the matching $\Quiv_\psi$ are exactly those cells corresponding to degree-maximal chains.
\end{lemma}
\begin{proof}
  First we will show that no cell appearing in an edge of $\Quiv_\psi$ can be degree-maximal and thus the cells associated to degree-maximal chains are critical cells. Suppose that $\alpha\leftarrow \beta$ is an edge in $\Quiv_\psi$ for  \[\alpha=\hull\{v_{\sigma_0},\ldots,\widehat{v}_{\sigma_s},\ldots,v_{\sigma_d}\}\text{ and } \beta=\hull\{v_{\sigma_0},\ldots,v_{\sigma_s},\ldots,v_{\sigma_{d}}\}.\]
  Then $\alpha$ cannot correspond to a degree-maximal chain, as $\sigma_{s-1}$ cannot be maximal among faces of $\sigma_{s+1}$ with weight $\psi(\sigma_{s-1})$. Moreover, $\beta$ is not degree-strict as $\psi(\sigma_{s-1}) = \psi(\sigma_s)$ and thus is also not degree-maximal.
  
  For the other direction, consider a chain $\sigma_0\prec\cdots\prec\sigma_d$ that is not degree-maximal. If $\sigma_d$ is not maximal among cells with weight $\psi(\sigma_d)$ then let $\sigma_{d}'$ be a cell maximal among cells with weight $\psi(\sigma_d)$. Then the pair
  \[\hull\{v_{\sigma_0},\ldots,v_{\sigma_d}\}\leftarrow \hull\{v_{\sigma_0},\ldots,v_{\sigma_{d}},v_{\sigma_{d}'}\}\]
  is an element of $\Quiv_\psi$.

  Otherwise, let $s>0$ be the minimal $s$ such that the chain $\sigma_s\prec \cdots \prec \sigma_d$ is degree-maximal. This implies that either $\psi(\sigma_{s-1})=\psi(\sigma_{s})$, or $\sigma_{s-1}$ is not maximal among faces of $\sigma_s$ with weight $\psi(\sigma_{s-1})$. In the first case, the edge \[\hull\{v_{\sigma_0},\ldots,\widehat{v}_{\sigma_s},\ldots,v_{\sigma_d}\}\leftarrow \hull\{v_{\sigma_0},\ldots,v_{\sigma_s},\ldots,v_{\sigma_{d}}\}\]
  is an element of $\Quiv_\psi$. In the second case, by Lemma~\ref{lem:unique-maximal-cell} there exists a face $\sigma_{s-1}'$ of $\sigma_s$ which is maximal among faces of $\sigma_s$ with weight $\psi(\sigma_{s-1})$. Then the edge
  \[\hull\{v_{\sigma_0},\ldots,v_{\sigma_d}\}\leftarrow \hull\{v_{\sigma_0},\ldots,v_{\sigma_{s-1}},v_{\sigma_{s-1}'},v_{\sigma_s},\ldots,v_{\sigma_{d}}\}\]
  is an edge in $\Quiv_\psi$.

  Thus the critical cells of the matching $\Quiv_\psi$ are exactly the cells of the barycentric subdivision $(\Lambda(L)_{\RR})^{BC}$ which correspond to degree-maximal chains.
\end{proof}

Having described the critical cells of the acyclic matching, denote by
\[
[\sigma_0\prec \cdots\prec \sigma_d]
\]
the cell of $\homotopy(\Quiv_\psi)(\Lambda(L)_{\RR}^{BC})$ corresponding to a critical cell $\hull\{v_{\sigma_0},\ldots,v_{\sigma_d}\}$ of the matching $\Quiv_\psi$. Now we describe the boundaries of these cells and connect these to the geometric realization of the dual complex of $\psi$.

\begin{lemma}
  \label{lem:unique-paths}
  Fix a shifted ceiling stratification $\psi$.
  For any pair of cells $\alpha_1,\alpha_m\subseteq \Lambda(L)_{\RR}^{BC}$ associated to degree-strict chains, there exists at most one path in $\graph^{\Quiv_\psi}(\Lambda(L)_{\RR}^{BC})$ of the form
  \[\alpha_1\rightarrow \beta_1 \rightarrow \alpha_2 \rightarrow \cdots \rightarrow \beta_{m-1} \rightarrow \alpha_m\]
  where each $\alpha_i\leftarrow \beta_i$ is an edge in $\Quiv_\psi$ for $i=1,\ldots,m-1$.

  Moreover, if such a path exists then the degree sequences associated to $\alpha_1$ and $\alpha_m$ must be identical.
\end{lemma}
\begin{proof}
  Observe in such a path, if $\alpha_i$ is not degree-strict then $\alpha_j$ will not be degree-strict for any $j>i$ as $\beta_i$ always has at least one more duplicated degree than $\alpha_i$ and so $\alpha_{i+1}$ cannot have fewer duplicated degrees than $\alpha_i$. As such, it suffices to consider paths where each $\alpha_i$ is degree-strict.

  Observe that $\Quiv_\psi$ is a matching, so for each $\alpha_i$ there is at most one $\alpha_i\leftarrow \beta_i$ in $\Quiv_\psi$, so it suffices to show that for each $\beta_i$ there is a unique edge $\beta_{i}\rightarrow \alpha_{i+1}$ such that $\alpha_{i+1}$ is degree-strict and $\alpha_{i+1}\leftarrow \beta_{i+1}$ is an edge in $\Quiv_\psi$ for some $\beta_{i+1}$.

  Let $\beta_i=\hull\{v_{\sigma_0},\ldots,v_{\sigma_d}\}$ indexed such that
  $\sigma_0\prec \cdots\prec \sigma_d$
  is a chain and let $s$ be such that the degree-maximal tail of that chain is
  $\sigma_{s}\prec \cdots \prec \sigma_d.$
  Then as $\alpha_i\leftarrow \beta_i$ is an edge in $\Quiv_\psi$, we have that $s>0$ and $\psi(\sigma_{s-1})=\psi(\sigma_s)$.
  But as $\alpha_{i+1}$ is degree-strict and cannot be equal to $\alpha_i=\hull\{v_{\sigma_0},\ldots,\widehat{v}_{\sigma_s},\ldots,v_{\sigma_d}\}$, this implies
  \[\alpha_{i+1}=\hull\{v_{\sigma_0},\ldots,\widehat{v}_{\sigma_{s-1}},\ldots,v_{\sigma_d}\}.\]
  Observe that the degree sequence associated to $\alpha_{i+1}$ is identical to the degree sequence associated to $\alpha_{i}$. Assuming the existence of the path described in the lemma, this implies that the degree sequences associated to $\alpha_1$ and $\alpha_m$ are identical.

  As this choice for $\alpha_{i+1}$ is unique, this implies that there exists at most one path of the form described.
\end{proof}

\begin{lemma}
  \label{lem:unique-degree-maximal}
  Given a degree sequence $\bm{a}_0<\cdots<\bm{a}_d$ with $\bm{a}_i\in \img \psi$, there exists at most one degree-maximal chain $\sigma_0\prec \cdots \prec \sigma_d$ satisfying $\psi(\sigma_i) = \bm{a}_i$.
\end{lemma}
\begin{proof}
  This is an immediate consequence of Lemma~\ref{lem:strata-maximal-cells} and Lemma~\ref{lem:unique-maximal-cell}.
\end{proof}

\begin{proposition}
  \label{prop:boundary-and-homotopy}
  The $d$-cells of $\homotopy(\Quiv_\psi)(\Lambda(L)_{\RR}^{BC})$ are in bijection with degree-maximal chains of cells. Moreover, for the cell $[\sigma_0\prec \cdots \prec \sigma_d]$, the cells appearing as facets in the boundary are exactly:
  \[
    [\sigma_0\prec \cdots \prec \sigma_{s-1}\prec \sigma_{s+1}\prec \cdots \prec \sigma_d]
    \]
    for $0\leq s < d$ or
    \[
    [\sigma_0'\prec \cdots \prec \sigma_{d-1}']
    \]
  where $\sigma_{d-1}'$ is the maximal cell with $\psi(\sigma_{d-1}')=\psi(\sigma_{d-1})$ and  $\sigma_{i}'$ is a cell maximal with $\psi(\sigma_{i}')=\psi(\sigma_{i})$ such that $\sigma_i'$ is a face of $\sigma_{i+1}'$ for $0\leq i <d-1$.

  In particular, there exists an $L$-invariant homotopy equivalence between $\Lambda(L)_{\RR}^{BC}$ and the geometric realization of the dual complex of $\psi$.
\end{proposition}
\begin{proof}
  The critical cells are as described in Lemma~\ref{lem:critical-cells}. Let $\tau = [\sigma_0\prec \cdots \prec \sigma_d]$ be a $d$-cell of $\homotopy(\Quiv_\psi)(\Lambda(L)_{\RR}^{BC})$. By Lemma~\ref{lem:degree-maximal-subchains}, for all $0\leq s< d$ the chain $\sigma_0\prec \cdots\prec \sigma_{s-1}\prec \sigma_{s+1}\prec \cdots \prec \sigma_d$ is degree-maximal and so the acyclic matching has no effect on this pair of cells and so the cell \[[\sigma_0\prec \cdots\prec \sigma_{s-1}\prec \sigma_{s+1}\prec \cdots \prec \sigma_d]\] appears in the boundary of $\tau$.

  Otherwise, as $\sigma_0\prec \cdots \prec \sigma_{d-1}$ is degree-strict, applying Lemma~\ref{lem:unique-paths} and Lemma~\ref{lem:unique-degree-maximal} gives the desired result. In particular the chain $\sigma_0'\prec \cdots \prec \sigma_{d-1}'$ described in the statement of the proposition is the unique degree-maximal chain with the same degree sequence as $\sigma_0\prec \cdots \prec \sigma_{d-1}$. Moreover, for $t=0,\ldots,d-1$, the following edge is in $\Quiv_\psi$:
  \[\hull \{ v_{\sigma_0}, \ldots,v_{\sigma_{t}},\ldots,v_{\sigma_{d-1}'}\} \leftarrow
  \hull \{ v_{\sigma_0}, \ldots,v_{\sigma_{t}},v_{\sigma_{t}'},\ldots,v_{\sigma_{d-1}'}\}.
  \]
  Consequently, applying Lemma~\ref{lem:unique-paths}, the only other cell appearing in the boundary of $\tau$ is $[\sigma_0'\prec \cdots \prec \sigma'_{d-1}]$.
  
  Additionally, the identification of degree-maximal chains with degree sequences in Lemma~\ref{lem:unique-degree-maximal} gives a bijection of cells between
  $\homotopy(\Quiv_\psi)(\Lambda(L)_{\RR}^{BC})$ and $\Dual_{\psi}$ compatible with the face poset.
  By Lemma~\ref{lem:unique-paths} the paths to the boundary facets through $\graph^{\Quiv_\psi}(\Lambda(L)_{\RR}^{BC})$ are unique,
  so the resulting CW complex $\homotopy(\Quiv_\psi)(\Lambda(L)_{\RR}^{BC})$ is a regular.
  It follows that $\homotopy(\Quiv_\psi)(\Lambda(L)_{\RR}^{BC})$ is a geometric realization of the dual complex of $\psi$.
  Finally, as the acyclic matching $\Quiv_\psi$ is $L$-invariant, the resulting homotopy equivalence will also be $L$-invariant.
\end{proof}

Now define a compatible $\ZZ^{2n}$-stratification on the dual complex by
\[\psi^{dual}(\{\bm{a}_0,\ldots,\bm{a}_d\}) = (\bm{a}_d,-\bm{a}_0).\]
The content of Proposition~\ref{prop:boundary-and-homotopy} then implies that this stratification is exactly the result of applying a homotopy arising from $\Quiv_\psi$ to the stratification $\psi^{BC}$.
Using this stratification, we conclude this section with the proof of Theorem~\ref{thm:dual-complex}.

\begin{proof}[Proof of Theorem~\ref{thm:dual-complex}]
  That the dual complex of $\psi$ is $L$-invariant is immediate as for a compatible $\ZZ^n$-stratification, the cell structure is $L$-invariant and the stratification function $\psi$ is $L$-equivariant.
  
  By Proposition~\ref{prop:boundary-and-homotopy}, there is an $L$-invariant homotopy between $\Lambda(L)_{\RR}^{BC}$ and the geometric realization of the dual complex. As $\Lambda(L)_{\RR}$ is homeomorphic to $L_{\RR}$ and $L_{\RR}/L$ is a torus, the quotient of the dual complex by $L$ has the homotopy type of a torus.
  
  That the complex $F_{\psi^{dual}}$ is a free resolution of $M_{\psi^{BC}}$ follows from observing that $\psi^{BC}(\alpha)=\psi^{BC}(\beta)$ for any edge $\alpha\leftarrow \beta$ in $ \Quiv_\psi$ and thus the acyclic matching $\Quiv_\psi$ is homogeneous. As $\psi^{dual}$ arises from $\psi^{BC}$ by applying the homotopy $\homotopy(\Quiv_\psi)$, the resolution $F_{\psi^{dual}}$ will be exactly the resolution that \cite{Batzies}*{Theorem~3.3.1} gives when applied to $F_{\psi^{BC}}$ viewed as a resolution supported on $\Lambda(L)_{\RR}^{BC}$. This implies that $\pi_{\Lambda(L)}(F_{\psi^{BC}})$ resolves $\Mepsilon = \pi_{\Lambda(L)}(M_{\psi^{BC}})$ as desired.
\end{proof}

\subsection*{The Case of the Second Hirzebruch Surface}

  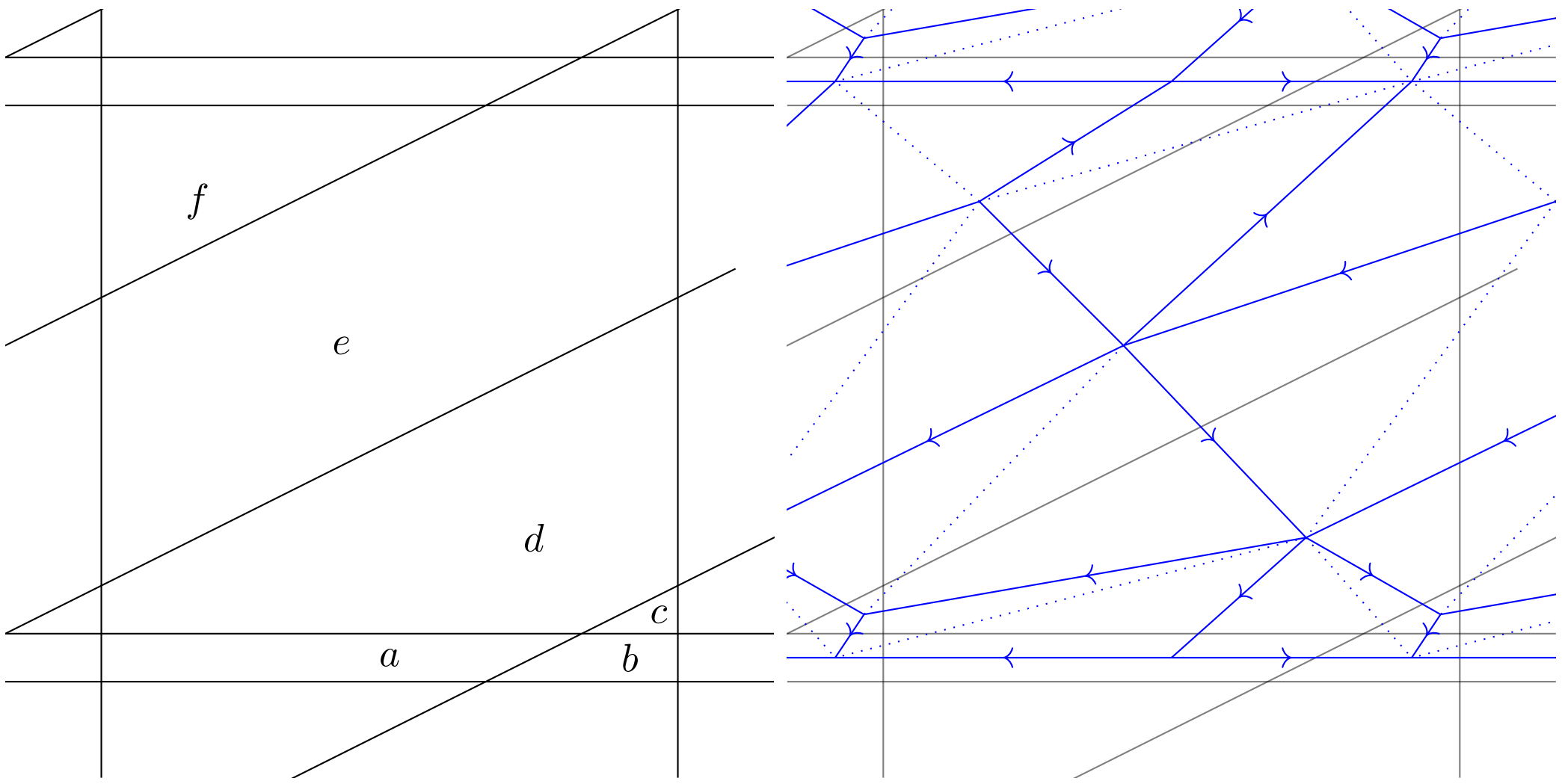
\begin{figure}[b]
    \begin{tikzpicture}[scale=0.85]
      \clip (-1,-1) rectangle (7,7);
      \coordinate (origin) at (0,0);
      \coordinate (mr) at (6,3);
      \coordinate (br) at (6,0);
      \coordinate (tr) at (6,6);
      \coordinate (tl) at (0,6);
      \coordinate (ml) at (0,3);

      \coordinate (a) at (3,0.25);
      \coordinate (b) at (5.5,0.25);
      \coordinate (c) at (5.8,0.7);
      \coordinate (d) at (4.5,1.5);
      \coordinate (e) at (2.5,3.5);
      \coordinate (f) at (1,5);

      \node at (a) {$a$};
      \node at (b) {$b$};
      \node at (c) {$c$};
      \node at (d) {$d$};
      \node at (e) {$e$};
      \node at (f) {$f$};
      
      \draw ($(origin)-(0,1)$) -- ($(tl)+(0,1)$);
      \draw ($(origin)-(1,0)$) -- ($(br)+(1,0)$);
      \draw ($(br)-(0,1)$) -- ($(tr)+(0,1)$);
      \draw ($(tl)-(1,0)$) -- ($(tr)+(1,0)$);

      \begin{scope}
        \draw ($(origin)-(1,0)+(0,0.5)$) -- ($(br)+(1,0)+(0,0.5)$);
        \draw ($(tl)-(1,0)+(0,0.5)$) -- ($(tr)+(1,0)+(0,0.5)$);
      \end{scope}
      
      \begin{scope}
        \draw ($(ml)-(1,0.5)+(-0.4,0.8)$) -- ($(tr)+(1,0.5)+(-0.4,0.8)$);
        \draw ($(origin)-(1,0.5)+(-0.4,0.8)$) -- ($(mr)+(1,0.5)+(-0.4,0.8)$);
        \draw ($(br)-(4,2)+(-0.4,0.8)$) -- ($(br)+(2,1)+(-0.4,0.8)$);
        \draw ($(tl)-(1,0.5)+(-0.4,0.8)$) -- ($(tl)+(1,0.5)+(-0.4,0.8)$);
      \end{scope}
    \end{tikzpicture}
    \begin{tikzpicture}[scale=0.85]
      \clip (-1,-1) rectangle (7,7);
      \begin{scope}[opacity=0.5]
        \coordinate (origin) at (0,0);
        \coordinate (mr) at (6,3);
        \coordinate (br) at (6,0);
        \coordinate (tr) at (6,6);
        \coordinate (tl) at (0,6);
        \coordinate (ml) at (0,3);

        \draw ($(origin)-(0,1)$) -- ($(tl)+(0,1)$);
        \draw ($(origin)-(1,0)$) -- ($(br)+(1,0)$);
        \draw ($(br)-(0,1)$) -- ($(tr)+(0,1)$);
        \draw ($(tl)-(1,0)$) -- ($(tr)+(1,0)$);

        \begin{scope}
          \draw ($(origin)-(1,0)+(0,0.5)$) -- ($(br)+(1,0)+(0,0.5)$);
          \draw ($(tl)-(1,0)+(0,0.5)$) -- ($(tr)+(1,0)+(0,0.5)$);
        \end{scope}
        
      \begin{scope}
        \draw ($(ml)-(1,0.5)+(-0.4,0.8)$) -- ($(tr)+(1,0.5)+(-0.4,0.8)$);
        \draw ($(origin)-(1,0.5)+(-0.4,0.8)$) -- ($(mr)+(1,0.5)+(-0.4,0.8)$);
        \draw ($(br)-(4,2)+(-0.4,0.8)$) -- ($(br)+(2,1)+(-0.4,0.8)$);
        \draw ($(tl)-(1,0.5)+(-0.4,0.8)$) -- ($(tl)+(1,0.5)+(-0.4,0.8)$);
      \end{scope}
      \end{scope}
      \begin{scope}[color = blue,decoration={
            markings,
        mark=at position 0.5 with {\arrow{>}}}]
        \coordinate (a) at (3,0.25);
        \coordinate (b) at (5.5,0.25);
        \coordinate (c) at (5.8,0.7);
        \coordinate (d) at (4.4,1.5);
        \coordinate (e) at (2.5,3.5);
        \coordinate (f) at (1,5);

        \draw[postaction={decorate}] (a) -- (b);
        \draw[postaction={decorate}] (a) -- ($(b)+(-6,0)$);
        \draw[postaction={decorate}] (c) -- (b);
        \draw[postaction={decorate}] (d) -- (c);
        \draw[postaction={decorate}] (d) -- (a);
        \draw[postaction={decorate}] (d) -- ($(c)+(-6,0)$);
        \draw[postaction={decorate}] (e) -- (d);
        \draw[postaction={decorate}] (e) -- ($(d)+(-6,0)$);
        \draw[postaction={decorate}] (e) -- ($(b)+(0,6)$);
        \draw[postaction={decorate}] (f) -- (e);
        \draw[postaction={decorate}] (f) -- ($(e)+(-6,0)$);
        \draw[postaction={decorate}] (f) -- ($(a)+(0,6)$);

        \draw[postaction={decorate}] ($(a)+(6,0)$) -- ($(b)+(6,0)$);
        \draw[postaction={decorate}] ($(a)+(6,0)$) -- ($(b)+(0,0)$);
        \draw[postaction={decorate}] ($(c)+(6,0)$) -- ($(b)+(6,0)$);
        \draw[postaction={decorate}] ($(d)+(6,0)$) -- ($(c)+(6,0)$);
        \draw[postaction={decorate}] ($(d)+(6,0)$) -- ($(a)+(6,0)$);
        \draw[postaction={decorate}] ($(d)+(6,0)$) -- ($(c)+(0,0)$);
        \draw[postaction={decorate}] ($(e)+(6,0)$) -- ($(d)+(6,0)$);
        \draw[postaction={decorate}] ($(e)+(6,0)$) -- ($(d)+(0,0)$);
        \draw[postaction={decorate}] ($(e)+(6,0)$) -- ($(b)+(6,6)$);
        \draw[postaction={decorate}] ($(f)+(6,0)$) -- ($(e)+(6,0)$);
        \draw[postaction={decorate}] ($(f)+(6,0)$) -- ($(e)+(0,0)$);
        \draw[postaction={decorate}] ($(f)+(6,0)$) -- ($(a)+(6,6)$);

        \draw[postaction={decorate}] ($(a)+(-6,0)$) -- ($(b)+(-6,0)$);
        \draw[postaction={decorate}] ($(a)+(-6,0)$) -- ($(b)+(-12,0)$);
        \draw[postaction={decorate}] ($(c)+(-6,0)$) -- ($(b)+(-6,0)$);
        \draw[postaction={decorate}] ($(d)+(-6,0)$) -- ($(c)+(-6,0)$);
        \draw[postaction={decorate}] ($(d)+(-6,0)$) -- ($(a)+(-6,0)$);
        \draw[postaction={decorate}] ($(d)+(-6,0)$) -- ($(c)+(-12,0)$);
        \draw[postaction={decorate}] ($(e)+(-6,0)$) -- ($(d)+(-6,0)$);
        \draw[postaction={decorate}] ($(e)+(-6,0)$) -- ($(d)+(-12,0)$);
        \draw[postaction={decorate}] ($(e)+(-6,0)$) -- ($(b)+(-6,6)$);
        \draw[postaction={decorate}] ($(f)+(-6,0)$) -- ($(e)+(-6,0)$);
        \draw[postaction={decorate}] ($(f)+(-6,0)$) -- ($(e)+(-12,0)$);
        \draw[postaction={decorate}] ($(f)+(-6,0)$) -- ($(a)+(-6,6)$);

        \draw[postaction={decorate}] ($(a)+(0,6)$) -- ($(b)+(0,6)$);
        \draw[postaction={decorate}] ($(a)+(0,6)$) -- ($(b)+(-6,6)$);
        \draw[postaction={decorate}] ($(c)+(0,6)$) -- ($(b)+(0,6)$);
        \draw[postaction={decorate}] ($(d)+(0,6)$) -- ($(c)+(0,6)$);
        \draw[postaction={decorate}] ($(d)+(0,6)$) -- ($(a)+(0,6)$);
        \draw[postaction={decorate}] ($(d)+(0,6)$) -- ($(c)+(-6,6)$);

        \draw[postaction={decorate}] ($(a)+(-6,6)$) -- ($(b)+(-6,6)$);
        \draw[postaction={decorate}] ($(c)+(-6,6)$) -- ($(b)+(-6,6)$);
        \draw[postaction={decorate}] ($(d)+(-6,6)$) -- ($(c)+(-6,6)$);
        \draw[postaction={decorate}] ($(d)+(-6,6)$) -- ($(a)+(-6,6)$);

        \draw[postaction={decorate}] ($(a)+(6,6)$) -- ($(b)+(6,6)$);
        \draw[postaction={decorate}] ($(a)+(6,6)$) -- ($(b)+(0,6)$);
        \draw[postaction={decorate}] ($(c)+(6,6)$) -- ($(b)+(6,6)$);
        \draw[postaction={decorate}] ($(d)+(6,6)$) -- ($(c)+(6,6)$);
        \draw[postaction={decorate}] ($(d)+(6,6)$) -- ($(a)+(6,6)$);
        \draw[postaction={decorate}] ($(d)+(6,6)$) -- ($(c)+(0,6)$);

        
        \draw[dotted] (d)--(b);
        \draw[dotted] (d)--($(b)+(-6,0)$);
        \draw[dotted] (e)--($(c)+(-6,0)$);
        \draw[dotted] (f)--($(d)+(-6,0)$);
        \draw[dotted] (f)--($(b)+(0,6)$);
        \draw[dotted] (f)--($(b)+(-6,6)$);
        \draw[dotted] ($(f)+(6,0)$)--($(d)+(0,0)$);
        \draw[dotted] ($(f)+(6,0)$)--($(b)+(0,6)$);
        \draw[dotted] ($(d)+(0,6)$)--($(b)+(-6,6)$);
        \draw[dotted] ($(d)+(-6,0)$)--($(b)+(-6,0)$);
        \draw[dotted] ($(d)+(0,6)$)--($(b)+(0,6)$);
        \draw[dotted] ($(d)+(6,6)$)--($(b)+(0,6)$);
        \draw[dotted] ($(d)+(6,0)$)--($(b)+(0,0)$);
        \draw[dotted] ($(e)+(0,6)$)--($(c)+(-6,6)$);
        \draw[dotted] ($(e)+(6,6)$)--($(c)+(0,6)$);
        \draw[dotted] ($(e)+(6,0)$)--($(c)+(0,0)$);
      \end{scope}
    \end{tikzpicture}
    \caption{An $\epsilon$ shifted stratification of the second Hirzebruch surface along the corresponding dual complex.}
    \label{fig:h2-dual-complex}
  \end{figure}

  For additional clarity, we work out the case of the second Hirzebruch surface in detail.
In this case $R = \field[x_1,x_2,x_3,x_4]$, $S=\field[x_1,\ldots,x_4,y_1,\ldots,y_4]$
with $\deg(x_1)=\deg(x_3)=(1,0)$, $\deg(x_2) = (-2,1)$, and $\deg(x_4)=(0,1)$.
With $\epsilon = (0,1/12,1/3,0)$, we get the stratification depicted on the left of Figure~\ref{fig:h2-dual-complex}.
On the right the quiver has been denoted by solid blue arrows superposed on the stratification,
and the additional edges of the dual complex are depicted as dotted blue edges.
Because the hyperplanes are transverse, the unique maximal cell of each stratum is a maximal cell of $L_{\RR}$,
labeled in the figure using $a$ through $f$.

Then the monomials arising as $x^{\psi_\epsilon(\sigma)}$ for each maximal cell are listed below:
  \[\begin{array}{c|cccccc}
    \text{Cell} & a & b & c & d & e & f\\\hline
    \text{Monomial} & x_1 & \frac{x_1}{x_3} & \frac{x_1x_2}{x_3} & x_1x_2 & x_1x_2x_3 & x_1x_2x_3^2\\
  \end{array}\]

  The chain complex that results from the dual complex construction is as follows:
  \begin{equation*}
    \begin{split}
      &
      S(-1,0,1,0)
      \oplus
      S(0,0,0,0)
      \oplus
      S(2,-1,-2,1)
      \oplus
      S(1,-1,-1,1)
      \oplus
      S(0,-1,0,1)
      \oplus
      S(-1,-1,1,1)
    \\[2pt]
    &\qquad \xleftarrow{\partial_1}
    \begin{matrix}
      S(-1,0,0,0)^2
      \oplus
      S(2,-1,0,0)
      \oplus
      S(1,-1,1,0) 
      \oplus
      S(1,-1,-2,1)^2 
      \oplus\\
      S(0,-1,0,0) 
      \oplus
      S(0,-1,-1,1)^2 
      \oplus
      S(-1,-1,1,0) 
      \oplus
      S(-1,-1,0,1)^2 
      \oplus\\
      S(1,-1,0,0)^2 
      \oplus
      S(0,-1,-2,1) 
      \oplus
      S(-1,-1,0,0)^2 
      \oplus
      S(-1,-1,-1,1) 
      \\
    \end{matrix}\\[2pt]
    & \qquad \xleftarrow{\partial_2}
    \begin{matrix}
      S(1,-1,0,0)^4 
      \oplus
      S(0,-1,-2,1)^2 
      \oplus
      S(-1,-1,0,0)^4 
      \oplus 
      S(-1,-1,-1,1)^2 
      \\
    \end{matrix} \leftarrow 0
    \end{split}
  \end{equation*}

  Then the matrices of the differentials (with an appropriate choice of orientation) are

  \[
  \partial_1=\begin{psmallmatrix}
    -y_1 & -y_3 & 0 & x_2 & 0 & 0 & 0 & 0 & 0 & x_4 & 0 & 0 & 0 & 0 & 0 & 0 & 0 & 0\\
    x_1 & x_3 & x_2 & 0 & 0 & 0 & x_4 & 0 & 0 & 0 & 0 & 0 & x_1x_2 & x_2x_3 & 0 & x_1x_4 & x_3x_4 & 0\\
    0 & 0 & -y_2 & 0 & x_1 & x_3 & 0 & 0 & 0 & 0 & 0 & 0 & 0 & 0 & x_1x_3 & 0 & 0 & 0\\
    0 & 0 & 0 & -y_2 & -y_1 & -y_3 & 0 & x_1 & x_3 & 0 & 0 & 0 & -y_1y_2 & -y_2y_3 & 0 & 0 & 0 & x_1x_3\\
    0 & 0 & 0 & 0 & 0 & 0 & -y_4 & -y_1 & -y_3 & 0 & x_1 & x_3 & 0 & 0 & -y_1y_3 & 0 & 0 & 0\\
    0 & 0 & 0 & 0 & 0 & 0 & 0 & 0 & 0 & -y_4 & -y_1 & -y_3 & 0 & 0 & 0 & -y_1y_4 & -y_3y_4 & -y_1y_3\\
  \end{psmallmatrix}
  \]
    \[
    \partial_2=\begin{psmallmatrix}
    x_2 & 0 & 0 & 0 & 0 & 0 & x_4 & 0 & 0 & 0 & 0 & 0\\ 
    0 & x_2 & 0 & 0 & 0 & 0 & 0 & x_4 & 0 & 0 & 0 & 0\\ 
    0 & 0 & x_1 & x_3 & 0 & 0 & 0 & 0 & 0 & 0 & 0 & 0\\ 
    y_1 & y_3 & 0 & 0 & 0 & 0 & 0 & 0 & 0 & 0 & 0 & 0\\ 
    0 & 0 & y_2 & 0 & x_3 & 0 & 0 & 0 & 0 & 0 & 0 & 0\\ 
    0 & 0 & 0 & y_2 & 0 & x_1 & 0 & 0 & 0 & 0 & 0 & 0\\ 
    0 & 0 & 0 & 0 & 0 & 0 & 0 & 0 & x_1 & x_3 & 0 & 0\\ 
    0 & 0 & 0 & 0 & 0 & y_3 & 0 & 0 & 0 & 0 & x_3 & 0\\ 
    0 & 0 & 0 & 0 & y_1 & 0 & 0 & 0 & 0 & 0 & 0 & x_1\\ 
    0 & 0 & 0 & 0 & 0 & 0 & y_1 & y_3 & 0 & 0 & 0 & 0\\ 
    0 & 0 & 0 & 0 & 0 & 0 & 0 & 0 & y_4 & 0 & 0 & y_3\\ 
    0 & 0 & 0 & 0 & 0 & 0 & 0 & 0 & 0 & y_4 & y_1 & 0 \\ 
    -1 & 0 & -1 & 0 & 0 & 0 & 0 & 0 & 0 & 0 & 0 & 0\\ 
    0 & -1 & 0 & -1 & 0 & 0 & 0 & 0 & 0 & 0 & 0 & 0\\ 
    0 & 0 & 0 & 0 & -1 & -1 & 0 & 0 & 0 & 0 & 0 & 0\\ 
    0 & 0 & 0 & 0 & 0 & 0 & -1 & 0 & -1 & 0 & 0 & 0\\ 
    0 & 0 & 0 & 0 & 0 & 0 & 0 & -1 & 0 & -1 & 0 & 0\\ 
    0 & 0 & 0 & 0 & 0 & 0 & 0 & 0 & 0 & 0 & -1 & -1\\ 
  \end{psmallmatrix}
  \]

  This complex is not minimal and even has quadratic terms.
  However, it is not hard to see that there is an acyclic matching that removes all of the dotted edges in Figure~\ref{fig:h2-dual-complex}.
  After that acyclic matching, the complex will be linear.
  Moreover, the minimized complex matches with the complex given in \cite{brown-sayrafi}*{Example~3.9}.
  This generalizes easily to the remainder of the Hirzebruch surfaces, and likely to all of the Picard rank 2 toric varieties.
  Some care has to be taken in the choice of $\epsilon$ to agree exactly with \cite{brown-sayrafi}.
  In the case of Hirzebruch surfaces, using $\epsilon$ where the corresponding $\QQ$-divisor is $\QQ$-ample suffices.

\section{Comparison to Homotopy Path Algebras}
\label{sec:hpa}

Finally, we connect the resolution arising from the barycentric subdivision and the cellular resolutions arising from homotopy path algebras as described by \cite{FH25}*{Section~6.1}.
As the language is sufficiently different and the construction of the complexes is important, we reimplement a couple constructions from \cite{FH25}*{Section~3.1 and Section~6.1} here to emphasize certain details we require.

Akin to \cite{FH25}*{Definition~3.4}, we will construct homotopy path algebras from the quivers $Q_\psi$ and $\widetilde{Q}_\psi$.
Denote by $\field \widetilde{Q}_\psi$ the path algebra over the quiver $\widetilde{Q}_\psi$. We take the convention that the product $pq$ of paths $p$ and $q$ in $\widetilde{Q}_\psi$ is the path consisting of $p$ followed by $q$.

From Definition~\ref{def:monomial-quiver}, to each edge of $\widetilde{Q}_\psi$ is associated a monomial in $R$. Now to each path, associate a monomial by taking the product of all of the edge monomials. Construct a $\field$-algebra $\widetilde{A}$ as a quotient of $\field \widetilde{Q}_\psi$ by identifying paths with identical monomials. By Lemma~\ref{lem:quiver-basics} the path labels in $\widetilde{Q}_\psi$ are uniquely determined by the start and end points, and so $\widetilde{A}$ can also be constructed by identifying all paths with the same start and end points.

These relations are $L$-invariant as the monomial labels on the edges of the quiver are $L$-invariant. Consequently, they descend to relations on the path algebra $\field Q_\psi$. Denote the resulting homotopy path algebra by $A$.

It is immediate that $A$ is a homotopy path algebra as the relations satisfy both the left and right cancellative conditions of \cite{FH25}*{Definition~3.1}.
For $\widetilde{A}$, some care is required as \cite{FH25}*{Definition~3.7} is only stated for finite quivers, and $\widetilde{Q}_\psi$ is infinite. However, relaxing the definition in this circumstance will cause no issues. In particular, the construction of the cell complex $X_{\widetilde{A}}$ will only require that the quiver has finitely many paths between any two vertices and that there is a global bound on the length of any path. With these constraints, the cell complex $X_{\widetilde{A}}$ will be a finite dimensional, locally finite CW complex. The construction of the cellular cosheaf and corresponding complex in \cite{FH25}*{Section~6} extends to this case without further constraints.

As a clarification of notation, a path in $A$ is an element of $A$ corresponding to a path in $Q_\psi$.
Moreover, we define a partial order on paths in $Q_\psi$ where $p\leq q$ if $p$ is a prefix of $q$,
or equivalently if $q=pr$ for some path $r$ in $Q_\psi$. 
This partial order induces a partial order on the paths of $A$ where $p\leq q$
if there exist representative paths of $Q_\psi$ $p'$ and $q'$ with $p'\leq q'$. We denote the poset of paths in $A$ by $\Path_A$. Similarly, we denote the poset of paths in $\widetilde{A}$ by $\Path_{\widetilde{A}}$.
As an abuse of notation for a path $p$ in $A$, we denote by $h(p)$ and $t(p)$ the head and tail of $p$ respectively. Similarly, we define a partial order on paths in $\widetilde{Q}_\psi$ and $\widetilde{A}$.

Moreover, for paths $p,q$ in $A$ (or $\widetilde{A}$) satisfying $p\leq q$ denote by $q/p$ the unique path in $A$ such that $q=p(q/p)$. The definition of the partial order on $A$ ensures that such a path exists, and uniqueness is immediate from the left cancellative property of homotopy path algebras.

In a similar way to \cite{FH25}*{Definition~3.7 and Proposition~3.8}, the cell complex $X_{\widetilde{A}}$ is constructed as a quotient of the order complex on the poset $\Path_{\widetilde{A}}$ denoted $\Order(\Path_{\widetilde{A}})$.

While $X_{A}$ as constructed in \cite{FH25}*{Definition~3.7} is only a $\Delta$-complex (a semi-simplicial complex), applying the definition to $\widetilde{A}$ to obtain $X_{\widetilde{A}}$ yields a simplicial complex. We now replicate that construction in the case of $\widetilde{A}$ to explicitly realize $X_{\widetilde{A}}$ as a subcomplex of the simplex over the vertices of $\widetilde{Q}_\psi$.

To any chain of paths $p_0<p_1<\cdots<p_d$ in $\widetilde{A}$, we get a sequence of vertices of the quiver given by the end points of each path.
This set $\{h(p_d),\ldots,h(p_0)\}$ is a simplex on the vertices of the quiver $\widetilde{Q}_\psi$.
We reverse the order when writing to match the order on $\ZZ^n$ where the vertices live.
Then $X_{\widetilde{A}}$ is the quotient of $\Order(\Path_{\widetilde{A}})$ where two chains of paths $p_0<\cdots<p_d$ and $q_0<\cdots < q_d$ are identified if $h(p_i) = h(q_i)$ for all $i=0,\ldots,d$. To see that this is a well-defined simplicial complex, observe that the resulting simplicial complex is isomorphic to
\[\left\{\{h(p_d),\ldots,h(p_0)\} \,|\, \text{ for all chains } p_0<p_1<\cdots<p_d\text{ for paths } p_i\in \Path_{\widetilde{A}}\right\}.\]

For a cell $\sigma$ in $X_{A}$ or $X_{\widetilde{A}}$ corresponding to the chain of paths $p_0<p_1<\cdots<p_d$, define the head and tail by $h(\sigma)=h(p_d)$ and $t(\sigma)=h(p_0)$. The choice to use $t(\sigma)=h(p_0)$ is intentional and ensures that the tail of the cell is independent of representative in $\Order(\Path_{\widetilde{A}})$. In particular, if $p_0$ is a zero-length path, then $t(\sigma)=h(p_0)=t(p_0)$.

Moreover, the quotient construction of $X_{\widetilde{A}}$ from $\Order(\Path_{\widetilde{A}})$ is compatible with the $L$-action, so $X_{\widetilde{A}}$ acquires an $L$-action. We now prove that the quotient of $X_{\widetilde{A}}$ by the $L$-action recovers the space $X_{A}$ as defined in \cite{FH25}*{Definition~3.7}.
\begin{lemma}
  \label{lem:XA-quotient}
  The space $X_A$ as defined in \cite{FH25}*{Definition~3.7 and Proposition~3.8} is the quotient of $X_{\widetilde{A}}$ by the $L$-action induced by the action of $L$ on the quiver $\widetilde{Q}_\psi$.
\end{lemma}
\begin{proof}

  As $A$ is the quotient of $\widetilde{A}$ under the $L$-action, and the quiver $Q_\psi$ is acyclic, $\Order(\Path_{A})$ is the quotient of $\Order(\Path_{\widetilde{A}})$ under the $L$-action. Thus there is a surjective map from $\Order(\Path_{A})$ to the quotient of $X_{\widetilde{A}}$ under the $L$-action compatible with the cell structure.
  Thus to show that $X_A$ is the quotient of $X_{\widetilde{A}}$ by the $L$-action it suffices to show that the maps from $\Order(\Path_A)$ to both $X_{A}$ and the quotient identify the same simplices.

  Let $p_0<\cdots<p_d$ and $q_0<\cdots<q_d$ be chains in $\Order(\Path_A)$.
  Then these chains represent the same cell of $X_A$ if and only if $p_i/p_{i-1} = q_i/q_{i-1}$ in $A$ for $i=1,\ldots,d$. This follows from the identification of simplices given at the end of the proof of \cite{FH25}*{Proposition~3.8} and the observation that $(p'/q)/(p/q) = p'/p$ for paths $q<p<p'$.\footnote{The end of the proof of \cite{FH25}*{Proposition~3.8} contains a typo, the second simplex in the identification should have $e_{h(p_0)}$ as the first term of the chain. In particular, $e_{h(p_0)}$ is a prefix of $p_1/p_0$ where as $e_{t(p_0)}$ is not.}

  Now choose lifts
  $\widetilde{p}_0<\cdots < \widetilde{p}_d$ and $\widetilde{q}_0<\cdots < \widetilde{q}_d$
  to chains of paths in $\widetilde{A}$.

  Then as $A$ is the quotient of $\widetilde{A}$ by an $L$-action, $p_i/p_{i-1} = q_i/q_{i-1}$ if and only if
  $\widetilde{p}_i/\widetilde{p}_{i-1}$ and $\widetilde{q}_i/\widetilde{q}_{i-1}$ differ by an $L$-translation.
  Moreover, as the head of $\widetilde{p}_{i}/\widetilde{p}_{i-1}$ is the tail of $\widetilde{p}_{i+1}/\widetilde{p}_i$ and the action of $L$ on vertices is a free action, there exists a unique element $\bm{v}\in L$ such that $\widetilde{p}_i/\widetilde{p}_{i-1}=\widetilde{q}_i/\widetilde{q}_{i-1} + \bm{v}$ for all $i=1,\ldots,d$.

  As $\widetilde{A}$ has at most one path between any pair of vertices and the head of $\widetilde{p}_{i}/\widetilde{p}_{i-1}$ is the tail of $\widetilde{p}_{i+1}/\widetilde{p}_i$, the sequence of paths $\widetilde{p}_{1}/\widetilde{p}_0,\ldots,\widetilde{p}_{i+1}/\widetilde{p}_{0}$ is uniquely determined by the sequence of vertices $h(\widetilde{p}_0),\ldots,h(\widetilde{p}_d)$ in $\widetilde{A}$. Thus the two chains of paths determine the same cell in $X_A$ if and only if the sequence of vertices $h(\widetilde{p}_0),\ldots,h(\widetilde{p}_d)$ and $h(\widetilde{q}_0),\ldots,h(\widetilde{q}_d)$ differ by an $L$-translation. That is to say, two chains of paths determine the same cell if and only if they determine the same cell in $X_{\widetilde{A}}$ up to the $L$-action. Thus $X_A$ is the quotient of $X_{\widetilde{A}}$ by the $L$-action.
\end{proof}

In the framing of \cite{FH25}*{Section~6.1} cellular resolutions arise as cellular cosheaves. For clarity, we give a brief description of the connection between the stratifications, cellular cosheaves, and resolutions. For a regular cell complex $\mathcal{X}$, a cellular cosheaf $\mathcal{F}$ of $A$-bimodules is a presheaf on the opposite category of the cellular poset arising from the cell complex $\mathcal{X}$. For a discussion on these definitions and particularly why no gluing condition is required here, see \cite{Curry}*{Section~4 and Theorem 4.2.10}.

In particular, a cellular cosheaf assigns to each cell $\sigma$ of a CW complex $\mathcal{X}$ an $A$-bimodule $\mathcal{F}(\sigma)$ and for each pair of cells $\sigma'$ and $\sigma$ in $\mathcal{X}$ with $\sigma'\prec \sigma$, a map $\mathcal{F}(\sigma)\rightarrow \mathcal{F}(\sigma')$ such that for $\sigma''\prec \sigma'\prec \sigma$ the composite map $\mathcal{F}(\sigma)\rightarrow \mathcal{F}(\sigma')\rightarrow \mathcal{F}(\sigma'')$ agrees with the map $\mathcal{F}(\sigma)\rightarrow \mathcal{F}(\sigma'')$.
For our purposes, we observe that from a compatible $\ZZ^n$-stratification $\psi$ on $L_{\RR}$, we get a cellular cosheaf of $R$-modules on the cell structure of $L_{\RR}$ given by $\sigma \mapsto R\cdot x^{\psi(\sigma)}$ where the maps are given by the canonical inclusions arising as submodules of $\field[\ZZ^n]$.

For $X_{A}$ some more care is required, as it is not in general regular. One approach would be to use the entrance path category as seen in \cite{Curry}*{Proposition 11.2.7}.
However, for our purposes, it suffices to take a quotient.
As $X_A$ arises as a quotient of a regular CW complex $X_{\widetilde{A}}$ by a free $L$-action, we can apply the same quotient to the cellular category. Denote this quotient by $\Cell(X_{A})$.

From this cosheaf, we can construct a complex of $A$-modules, as seen in \cite{FH25}*{Equation~20}. By choosing an orientation on each cell of $X_A$ or $X_{\widetilde{A}}$, to each morphism $\partial$ in the cellular category, we can associate a sign $s(\partial)\in \{-1,0,1\}$. From an arbitrary cellular cosheaf and a choice of orientations, we get a complex of modules $F_\bullet$ with

\[F_i = \bigoplus_{\dim \sigma = i} \mathcal{F}(\sigma).\]
\[\partial_i(\sigma) = \sum_{\dim \tau = \dim \sigma - 1}\sum_{\partial \in \Hom(\tau,\sigma)} s(\partial) \mathcal{F}(\partial).\]
For any particular choice of orientation of the cells, the signs are chosen such that the above construction applied to the constant cosheaf $\underline{\field}$ recovers exactly the cellular chain complex.

Now we move to the definition of the cellular resolution of $\widetilde{A}$ as an $\widetilde{A}$-bimodule on $X_{\widetilde{A}}$. We will denote the resolution by $\widetilde{F}^{FH}$ and when necessary, denote the cosheaf by $\widetilde{\mathcal{F}}^{FH}$. Mirroring the notation of \cite{FH25}, for a vertex of $\widetilde{Q}_\psi$, define the following left and right projective $\widetilde{A}$-modules respectively:
\[\widetilde{P}_{\bm{v}} := \widetilde{A}e_{\bm{v}}\text{ and } \widetilde{P}_{\bm{v}}^{op} := e_{\bm{v}} \widetilde{A}.\]
Then, to a cell $\eta = [p_0<\cdots < p_d]\in X_{\widetilde{A}}$ assign the projective module
\[\widetilde{P}_{\eta} = \widetilde{P}_{h(p_0)}\boxtimes_{\field} \widetilde{P}^{op}_{h(p_d)}.\]
To define the cellular resolution, we also need to provide maps $\widetilde{P}_{\eta}\rightarrow \widetilde{P}_{\eta'}$ for each pair of cells $\eta,\eta'\in X_{\widetilde{A}}$ with $\eta'$ a facet of $\eta$. Again, we can mirror the setup of \cite{FH25}*{Equation~23} to define a map $\widetilde{P}_{\eta}\rightarrow \widetilde{P}_{\eta'}$ as follows.
\begin{equation}
  \label{eq:hpa-boundary}
x\otimes y\mapsto
\begin{cases}
  xp_1\otimes y & \text{if }\eta'=[p_1<\cdots<p_d]\\ 
  x\otimes y & \text{if }\eta'=[p_0<\cdots \widehat{p}_i \cdots < p_d] \text{ for } 0<i<d\\
  x\otimes (p_{d}/p_{d-1})y & \text{if }\eta'=[p_0<\cdots<p_{d-1}].
\end{cases}
\end{equation}

Again this cellular cosheaf is $L$-equivariant, so the quotient by $L$ gives a cosheaf on $X_A$ and thus a cellular resolution.
This cellular resolution exactly reproduces the resolution given by \cite{FH25}*{Section~6.1}.
Define the quotient cosheaf of $\widetilde{\mathcal{F}}^{FH}$, denoted $\mathcal{F}^{FH}$, as a cosheaf on $X_{A}$ as follows. On cells $\sigma$ of $X_A$ choose $\widetilde{\sigma}$ a lift to $X_{\widetilde{A}}$ then set $\mathcal{F}^{FH}(\sigma)=\widetilde{\mathcal{F}}^{FH}(\widetilde{\sigma})$. For morphisms $\gamma:\sigma\rightarrow \tau$ in $\Cell(X_A)$, fix a lift $\widetilde{\sigma}$ of $\sigma$. Then there exists a unique lift $\widetilde{\tau}$ such that the map lifts to a map $\widetilde{\gamma}:\widetilde{\sigma}\rightarrow \widetilde{\tau}$. Define the map $\mathcal{F}^{FH}(\gamma)$ to be the map $\widetilde{\mathcal{F}}^{FH}(\widetilde{\gamma})$.
\begin{lemma}
  \label{lem:FH-cosheaf-same}
  The cosheaf $\mathcal{F}^{FH}$ is a well-defined cosheaf of $A$-bimodules and is exactly the cosheaf defined in \cite{FH25}*{Section~6.1, Equations~21 and 22}.
\end{lemma}
\begin{proof}
  To see that it is a cosheaf of $A$-bimodules, observe that $\widetilde{A}$ has an $A$-bimodule structure.
  Viewing an element $a\in A$ as an equivalence class of elements of $\widetilde{A}$,
  define a left $A$-module structure on $\widetilde{A}$ by $a\cdot x = \sum_{b\in a} b\cdot x$ for $x\in \widetilde{A}$.
  This sum is finite since there are a finite number of elements of $\widetilde{A}$ where $b\cdot x \neq 0$.
  The right $A$-module structure is defined similarly and the compatibility of the left and right module structure is inherited from the $\widetilde{A}$-bimodule structure.
  
  Fix a cell $\sigma$ in $X_A$ and a lift $\widetilde{\sigma}$ to $X_{\widetilde{A}}$. We will show that as $A$-bimodules, $\widetilde{P}_{\widetilde{\sigma}} = P_\sigma$.
  Recall $P_{\sigma} =P_{t}\boxtimes_{\field} P^{op}_{h}$, where $h$ and $t$ are the head and tail respectively of the path corresponding to $\sigma$.
  Similarly $\widetilde{P}_{\widetilde{\sigma}} = \widetilde{P}_{\widetilde{t}}\boxtimes_{\field} \widetilde{P}^{op}_{\widetilde{h}}$, where $\widetilde{h}$ and $\widetilde{t}$ are the head and tail respectively of the path corresponding to $\widetilde{\sigma}$.

  We start by showing $\widetilde{P}_{\widetilde{t}}=P_{t}$ as a left $A$-module.
  Observe that $\widetilde{P}_{\widetilde{t}}$ and $P_t$ are submodules of $\widetilde{A}$ and $A$ respectively,
  and the quotient map $\widetilde{A}\rightarrow A$ sends the generator of $\widetilde{P}_{\widetilde{t}}$ to the generator of $P_t$.
  Thus there is a canonical surjective map of left $A$-modules $\widetilde{P}_{\widetilde{t}}\rightarrow P_t$ induced by the quotient map $\widetilde{A}\rightarrow A$.
  Moreover, as every path in $A$ ending at $t$ lifts uniquely to a path in $\widetilde{A}$ ending at $\widetilde{t}$,
  the map is also injective and thus gives a canonical bijection between the two.
  The $P^{op}$ case is similar, and thus $\widetilde{P}_{\widetilde{\sigma}} = P_\sigma$ as desired.

  Moreover, the above argument implies that for any element $\bm{v}\in L$, there is a canonical isomorphism of $A$-bimodules from $\widetilde{P}_\sigma$ to $\widetilde{P}_{\sigma+\bm{v}}$.
  Moreover, the maps defining the cosheaf as given in \eqref{eq:hpa-boundary} are also $L$-equivariant, as for two cells $\eta$ and $\eta'$, $L$-translates of each other, the maps differ exactly by the $L$-action on $\widetilde{A}$. Thus the cosheaf $\widetilde{\mathcal{F}}^{FH}$ is $L$-equivariant and the cosheaf $\mathcal{F}^{FH}$ is well-defined.

  We have shown that the quotient is well-defined and that the cosheaf agrees on cells. We now only need to show that it agrees on maps. For this, observe that all the maps defined by \eqref{eq:hpa-boundary} reduce to left and right multiplication by paths in $\widetilde{A}$, which under the quotient become left and right multiplication by paths in $A$, which we can verify matches with \cite{FH25}*{Equation~22}.
\end{proof}

These two lemmas above will allow us to work primarily with $\widetilde{A}$, which will significantly simplify the arguments.

\subsection*{\texorpdfstring{$A$}{A}-modules to and from \texorpdfstring{$S$}{S}-modules}

To complete the proof of Theorem~\ref{thm:hpa-pdim},
we discuss the connection between the $S$-module picture and the $A$-module picture.
There are two key pieces which must be handled to give the correspondence between the construction in Section~\ref{sec:dual} and the construction in \cite{FH25}*{Section~6.1}.
First we must show that the projectives from \cite{FH25} tensor to free modules with an appropriate choice of $S$-module with an $A$-module structure.
Second we must show that the diagonal module we resolved is exactly the tensor of $\widetilde{A}$ against the module $\widetilde{T}\boxtimes_{\field} \widetilde{T}'$ which we will define below.

As has been the convention, let $R$ be the Cox ring for $X$. Then let $\widetilde{T}$ be the free $R$-module generated in degrees given by the negatives of each of the vertices of $\widetilde{Q}_\psi$. Specifically 
\begin{equation}
  \widetilde{T} := \bigoplus_{\bm{v}\in \img \psi} R\cdot x^{-\bm{v}} \cong \bigoplus_{\bm{v}\in \img \psi} R(\deg (\bm{v})).
\end{equation}
For convenience denote the basis element of $\widetilde{T}$ corresponding to $\bm{v}$ by $f_{\bm{v}}$. Also define $\widetilde{T}'$ as the following explicit presentation of the dual with finitely generated support.
\begin{equation}
  \widetilde{T}' := \bigoplus_{\bm{v}\in \img \psi} R\cdot x^{\bm{v}} \cong \bigoplus_{\bm{v}\in \img \psi} R(-\deg (\bm{v})).
\end{equation}
Denote the dual basis element corresponding to $\bm{v}$ by $f'_{\bm{v}}$. This is the dual with finitely generated support in that there is a map $\widetilde{T}\otimes_{R} \widetilde{T}'\rightarrow R$ given on pure tensors by
\[f_{\bm{v}}\otimes f_{\bm{w}}' \mapsto 
\begin{cases}
  1 & \text{if }\bm{v}=\bm{w}\\
  0 & \text{otherwise}.
\end{cases}
\]
This pairing induces an injective map $\widetilde{T}'\hookrightarrow \Hom(\widetilde{T},R)$ and if we define $\Hom_c(M,N)$ to be those maps from $M$ with finitely generated support, or more formally
\[\Hom_c(M,N) := \varinjlim_{M'\subseteq M, M'\text{ f.g.}d}\Hom(M',N),\]
then $\widetilde{T}'\cong \Hom_c(\widetilde{T},R)$.

\begin{lemma}
  \label{lem:T-end}
  Let $\widetilde{A}$ and $\widetilde{T}$ be defined as above, then the map that identifies the unique path $\bm{v}\rightarrow \bm{w}$ with the endomorphism \[
\mu_{\bm{v},\bm{w}}(f_{\bm{u}}) =
\begin{cases}
  x^{\bm{v}-\bm{w}}f_{\bm{v}} & \bm{u} = \bm{w}\\
  0 & \text{otherwise.}
\end{cases}
\]

gives an isomorphism $\widetilde{A}\cong \End(\widetilde{T})_{\bm{0}}$.
\end{lemma}
\begin{proof}
Observe that the $\mu_{\bm{v},\bm{w}}$ form a $\field$-basis of $\End(\widetilde{T})_{\bm{0}}$. Consider the map of $\field$-algebras $\widetilde{A}\rightarrow \End(\widetilde{T})_{\bm{0}}$ which for the element of $\widetilde{A}$ corresponding to a path $p$ in $\widetilde{Q}_\psi$ assigns the map $\mu_{t(p),h(p)}$. As the paths in $\widetilde{A}$ are uniquely determined by their start and end vertices, this map defines an isomorphism of $\field$-vector spaces. Moreover, it is immediate that the map is compatible with the composition operations and so is an isomorphism of $\field$-algebras between $\End(\widetilde{T})_{\bm{0}}$ and $\widetilde{A}$.
\end{proof}

The identification of $\widetilde{A}$ with $\End(\widetilde{T})_{\bm{0}}$ gives $\widetilde{T}$ a left $\widetilde{A}$-module structure. Similarly, $\widetilde{T}'$ has a right $\widetilde{A}$-module structure via the identification with a submodule of $\Hom(\widetilde{T},R)$.

Define $T=\widetilde{T}\otimes_{R[L]} R$ and $T'=\widetilde{T}'\otimes_{R[L]} R$. Then $T$ has a $\Cl(X)$-graded structure as an $R$-module as $\ZZ^n/L=\ZZ^n/M_X\cong \Cl(X)$ and the vertices of $\psi$ are given by the homomorphic image of $\img \psi$ in $\ZZ^n/L$. In this graded structure, $T$ is
\[T \cong \bigoplus_{\alpha \in \img\psi/L} R(\alpha).\]
And $T'$ is
\[T' \cong \bigoplus_{\alpha \in \img\psi/L} R(-\alpha).\]

Then $T\boxtimes_{\field} T'$ has a $\Cl(X)\times\Cl(X)$-graded $S$-module structure via the isomorphism $S\cong R\otimes_{\field} R$. Moreover, as $A$ and $T$ are quotients of $\widetilde{A}$ and $\widetilde{T}$ respectively, a similar argument to Lemma~\ref{lem:T-end} gives $A\cong \End(T)_{\bm{0}}$. Consequently, $T$ and $T'$ have a left and right $A$-module structure respectively. From these $A$-module structures, $T\boxtimes_{\field} T'$ acquires an $A$-bimodule structure. Additionally, this identification of $A$ with $\End(T)_{\bm{0}}$ connects this construction to that of \cite{FH25}*{Definition 5.3} as the endomorphisms of the corresponding direct sum of line bundles arise exactly as degree $\bm{0}$ endomorphisms of the $S$-module $T$.

Using notation similar to that of \cite{associative-algebras-book}*{Section~10.1},
set $A^{e}:= A\otimes_{\field} A^{op}$. Then an $A$-bimodule is a module with a left (or right) $A^{e}$-module structure.
Similarly, define $\widetilde{A}^e:= \widetilde{A}\otimes_{\field}\widetilde{A}^{op}$.

There are two functors which play an important role through the rest of this section: tensor over $A^{e}$ against $T\boxtimes_{\field} T'$ and tensor over $\widetilde{A}^e$ against $\widetilde{T}\boxtimes_{\field} \widetilde{T}'$.
In both cases, the resulting object does not have any left or right $A$-module structure. As $\widetilde{T}\boxtimes_{\field} \widetilde{T}'$ and $T\boxtimes_{\field} T'$ both have $S$-module structures, the resulting object will have an $S$-module structure. 

Then a key property of these tensors is that they send the projective modules $P_\sigma$ to appropriately graded rank 1 free $S$-modules.
\begin{lemma}
  \label{lem:P-tensor-tilting}
  Let $\sigma$ be a cell in $X_A$ and $\sigma'$ a representative of $\sigma$ in $X_{\widetilde{A}}$. 
  \[\widetilde{P}_{\sigma'}\otimes_{\widetilde{A}^e} (\widetilde{T}\boxtimes_{\field} \widetilde{T}') = S\cdot x^{t(\sigma')}y^{-h(\sigma')},\]
and
  \[P_\sigma\otimes_{A^e} (T\boxtimes_{\field} T') = S(-\deg (t(\sigma')),\deg (h(\sigma'))).\]
  Moreover, the map
  \[\widetilde{P}_{\sigma'}\otimes_{\widetilde{A}^e} (\widetilde{T}\boxtimes_{\field} \widetilde{T}')\rightarrow P_\sigma\otimes_{A^e} (T\boxtimes_{\field} T'),\]
  induced by the maps $\widetilde{P}_{\sigma'}\rightarrow P_\sigma$ and $(\widetilde{T}\boxtimes_{\field} \widetilde{T}') \rightarrow (T\boxtimes_{\field} T')$ is the isomorphism that coarsens the grading.
\end{lemma}
\begin{proof}
  Let $\bm{u}$ and $\bm{v}$ in $\ZZ^n$ be the tail and head of $\sigma'$ respectively and let $\alpha$ and $\beta$ be their respective classes in $\ZZ^n/L\cong \Cl(X)$.
  
  Observe that the left $\widetilde{A}$-module $\widetilde{P}_{\bm{u}}$ satisfies
  \[ \widetilde{T}' \otimes_{\widetilde{A}} \widetilde{P}_{\bm{u}} = \widetilde{T}' \otimes_{\widetilde{A}} \widetilde{A}e_{\bm{u}} = R\cdot x^{\bm{u}}.\]
  Similarly
  \[ \widetilde{P}^{op}_{\bm{v}} \otimes_{\widetilde{A}} \widetilde{T} = e_{\bm{v}}\widetilde{A} \otimes_{\widetilde{A}} \widetilde{T}  = R\cdot x^{-\bm{v}}.\]
  Then under the identification of $S$ with $R\otimes_{\field} R$,
  \[\widetilde{P}_{\sigma'}\otimes_{\widetilde{A}^{e}} (\widetilde{T}\boxtimes_{\field} \widetilde{T}') = S\cdot x^{\bm{u}}y^{-\bm{v}}.\]

  Similarly, as $\alpha$ and $\beta$ are the vertices of $Q_\psi$ corresponding to $\bm{u}$ and $\bm{v}$ viewed as an element of $\Cl(X)$ respectively, the left $A$-module $P_\alpha$ satisfies
  \[T'\otimes_A P_{\alpha} = T'\otimes_A Ae_{\alpha} \cong R(-\alpha)\]
  and the right $A$-modules $P^{op}_{\beta}$ satisfies
  \[P^{op}_{\beta}\otimes_A T = e_{\beta}A\otimes_A T \cong R(\beta).\]
  Consequently
  \[P_{\sigma}\otimes_{A^e} (T\boxtimes_{\field} T') = S(-\alpha,\beta).\]
  
  For the map, observe that the generator of
  $\widetilde{P}_\sigma\otimes_{\widetilde{A}^{e}} (\widetilde{T}\boxtimes_{\field} \widetilde{T}')$
  as an $S$-module is the pure tensor $(e_{\bm{u}}\otimes e_{\bm{v}}) \otimes (f_{\bm{u}} \otimes f_{\bm{v}}')$.
  Under the map $\widetilde{P}_\sigma\otimes_{\widetilde{A}^e} (\widetilde{T}\boxtimes_{\field} \widetilde{T}')\rightarrow P_\sigma\otimes_{A^e} (T\boxtimes_{\field} T')$
  this maps to the corresponding pure tensor
  $(e_\alpha \otimes e_\beta) \otimes (f_\alpha \otimes f_\beta')$
  which generates $P_\sigma\otimes_{A^e} (T\boxtimes_{\field} T')$.
  Consequently the map is an isomorphism of $S$-modules. That it coarsens the grading is then immediate.
\end{proof}

We need one last lemma to reduce to working with $\widetilde{F}^{FH}$.
\begin{lemma}
  \label{lem:FHquotient}
  \[\pi_{\Lambda(L)}(\widetilde{F}^{FH}\otimes_{\widetilde{A}^{e}} (\widetilde{T}\boxtimes_{\field} \widetilde{T}')) \cong F^{FH}\otimes_{A^e} (T\boxtimes_{\field} T').\]
\end{lemma}
\begin{proof}

  It suffices to show that the cellular cosheaves $\widetilde{\mathcal{F}}^{FH}\otimes_{\widetilde{A}^e} (\widetilde{T}\boxtimes_{\field} \widetilde{T}')$ and $\mathcal{F}^{FH}\otimes_{A^e} (T\boxtimes_{\field} T')$ agree. 
  
  The cellular cosheaf $\widetilde{\mathcal{F}}^{FH}\otimes_{\widetilde{A}^e} (\widetilde{T}\boxtimes_{\field} \widetilde{T}')$ is an $L$-equivariant cellular cosheaf of $S$-modules on $X_{\widetilde{A}}$, and this naturally gives a cellular cosheaf of $S[\Lambda(L)]$-modules on $X_{A}$. In particular, we can define the cosheaf on $X_A$ which for a cell $\sigma$ with representative $\sigma'$ in $X_{\widetilde{A}}$ takes the value
  \begin{align*}
    \bigoplus_{\bm{v}\in L}\widetilde{P}_{\sigma'+\bm{v}}\otimes_{\widetilde{A}^e} (\widetilde{T}\boxtimes_{\field} \widetilde{T}') &= \bigoplus_{\bm{v}\in L}S\cdot x^{t(\sigma'+\bm{v})}y^{-h(\sigma'+\bm{v})}\\
      &= \bigoplus_{\bm{v}\in L}S\cdot x^{t(\sigma')}y^{-h(\sigma')}x^{\bm{v}}y^{-\bm{v}}.
  \end{align*}
  Similarly, the maps between these $S[\Lambda(L)]$-modules are defined by taking the direct sum across all $L$-translates.

  Then as $(\bm{v},-\bm{v})\in \Lambda(L)$ the functor $\pi_{\Lambda(L)}$ sends this module to $S(-\deg(t(\sigma')),\deg(h(\sigma')))$. Thus the cosheaf $\pi_{\Lambda(L)}(\widetilde{\mathcal{F}}^{FH}\otimes_{\widetilde{A}^e} (\widetilde{T}\boxtimes_{\field} \widetilde{T}'))$ assigns $S$-modules to cells by
  \[\sigma\mapsto S(-\deg(t(\sigma')),\deg(h(\sigma'))).\]
  On the other hand $P_{\sigma}\otimes_{A^e} (T\boxtimes_{\field}T')=S(-\deg(t(\sigma')),\deg(h(\sigma')))$ and so the two cosheaves assign isomorphic modules to each cell.

  We now show that the maps agree. Let $\sigma$ and $\tau$ be cells of $X_A$ and fix a morphism in $\Cell(X_A)^{op}$ $\sigma\rightarrow \tau$. Then there exist $\sigma'$ and $\tau'$ representatives of $\sigma$ and $\tau$ respectively such that $\tau'$ is a facet of $\sigma'$. Recall that the map $\widetilde{P}_{\sigma'}\rightarrow \widetilde{P}_{\tau'}$ is multiplication (either on the left or the right) by a path in $\widetilde{Q}_\psi$ and that the map $P_\sigma\rightarrow P_{\tau}$ is multiplication by the corresponding path in $Q_\psi$, either on the left or the right.

  From this we have the following commutative diagram of $A$-bimodules.
  \[\begin{tikzcd}
    \widetilde{P}_{\sigma'} \arrow[r] \arrow[d] & \widetilde{P}_{\tau'} \arrow[d]\\
    P_\sigma\arrow[r] & P_{\tau}
  \end{tikzcd}\]

  Tensoring with the map $\widetilde{T}\boxtimes_{\field}\widetilde{T}'\rightarrow T\boxtimes_{\field} T'$ gives the following commutative diagram.
  \[\begin{tikzcd}
    \widetilde{P}_{\sigma'} \otimes_{\widetilde{A}^e} (\widetilde{T}\boxtimes_{\field}\widetilde{T}')  \arrow[r] \arrow[d] & \widetilde{P}_{\tau'} \otimes_{\widetilde{A}^e} (\widetilde{T}\boxtimes_{\field}\widetilde{T}')\arrow[d]\\
    P_\sigma \otimes_{A^e} (T\boxtimes_{\field} T') \arrow[r] & P_{\tau} \otimes_{A^e} (T\boxtimes_{\field} T')
  \end{tikzcd}\]

  As previously noted, the vertical maps are isomorphisms of $S$-modules that coarsen the grading.

  Extending the commutative diagram gives
  \[\begin{tikzcd}
  \pi_{\Lambda(L)}(\bigoplus_{\widetilde{\sigma}\in \sigma} \widetilde{P}_{\widetilde{\sigma}} \otimes_{\widetilde{A}^e} (\widetilde{T}\boxtimes_{\field}\widetilde{T}'))  \arrow[r] & \pi_{\Lambda(L)}(\bigoplus_{\widetilde{\tau}\in \tau} \widetilde{P}_{\widetilde{\tau}} \otimes_{\widetilde{A}^e} (\widetilde{T}\boxtimes_{\field}\widetilde{T}'))\\
  \bigoplus_{\widetilde{\sigma}\in \sigma} \widetilde{P}_{\widetilde{\sigma}} \otimes_{\widetilde{A}^e} (\widetilde{T}\boxtimes_{\field}\widetilde{T}')  \arrow[r] \arrow[u] & \bigoplus_{\widetilde{\tau}\in \tau} \widetilde{P}_{\widetilde{\tau}} \otimes_{\widetilde{A}^e} (\widetilde{T}\boxtimes_{\field}\widetilde{T}') \arrow[u]\\
  \widetilde{P}_{\sigma'} \otimes_{\widetilde{A}^e} (\widetilde{T}\boxtimes_{\field}\widetilde{T}')  \arrow[r] \arrow[d] \arrow[u] & \widetilde{P}_{\tau'} \otimes_{\widetilde{A}^e} (\widetilde{T}\boxtimes_{\field}\widetilde{T}')\arrow[d] \arrow[u]\\
  P_\sigma \otimes_{A^e} (T\boxtimes_{\field} T') \arrow[r] & P_{\tau} \otimes_{A^e} (T\boxtimes_{\field} T')
  \end{tikzcd}\]

  The composite maps $ \widetilde{P}_{\sigma'} \otimes_{\widetilde{A}^e} (\widetilde{T}\boxtimes_{\field}\widetilde{T}')\rightarrow \pi_{\Lambda(L)}(\bigoplus_{\widetilde{\sigma}\in \sigma} \widetilde{P}_{\widetilde{\sigma}} \otimes_{\widetilde{A}^e} (\widetilde{T}\boxtimes_{\field}\widetilde{T}'))$ are isomorphisms compatible with the boundary maps, and so the following diagram commutes.

  \[\begin{tikzcd}
  \pi_{\Lambda(L)}(\bigoplus_{\widetilde{\sigma}\in \sigma} \widetilde{P}_{\widetilde{\sigma}} \otimes_{\widetilde{A}^e} (\widetilde{T}\boxtimes_{\field}\widetilde{T}'))  \arrow[r] & \pi_{\Lambda(L)}(\bigoplus_{\widetilde{\tau}\in \tau} \widetilde{P}_{\widetilde{\tau}} \otimes_{\widetilde{A}^e} (\widetilde{T}\boxtimes_{\field}\widetilde{T}'))\\
  \widetilde{P}_{\sigma'} \otimes_{\widetilde{A}^e} (\widetilde{T}\boxtimes_{\field}\widetilde{T}')  \arrow[r] \arrow[d] \arrow[u,"\cong"] & \widetilde{P}_{\tau'} \otimes_{\widetilde{A}^e} (\widetilde{T}\boxtimes_{\field}\widetilde{T}')\arrow[d] \arrow[u,"\cong"]\\
  P_\sigma \otimes_{A^e} (T\boxtimes_{\field} T') \arrow[r] & P_{\tau} \otimes_{A^e} (T\boxtimes_{\field} T')
  \end{tikzcd}\]

  Thus we have our isomorphism \[\pi_{\Lambda(L)}(\widetilde{F}^{FH}\otimes_{\widetilde{A}^{e}} (\widetilde{T}\boxtimes_{\field} \widetilde{T}')) \cong F^{FH}\otimes_{A^e} (T\boxtimes_{\field} T').\]
\end{proof}

Next, we relate the dual complex of $\psi$ to $X_A$.

\begin{lemma}
  \label{lem:inclusion}
  There exists an injective map of simplicial complexes from the dual complex to $X_{\widetilde{A}}$ denoted $\iota$ such that, up to a compatible choice of orientation, the induced map on complexes $F_{\psi^{dual}}\rightarrow F^{FH}\otimes_{\widetilde{A}^e} (\widetilde{T}\boxtimes_{\field} \widetilde{T}')$ given by $\sigma \mapsto e_{\iota(\sigma)}\otimes (f_{h(\iota(\sigma))}\boxtimes f_{t(\iota(\sigma))}')$ is an inclusion.
\end{lemma}
\begin{proof}
  Let $\{\bm{a}_0,\ldots,\bm{a}_d\}$ be a simplex in the dual complex for $\psi$ with $\bm{a}_0<\cdots<\bm{a}_d$.
  Then define the path $p_i$ to be the unique path in $\widetilde{A}$ from $\bm{a}_d$ to $\bm{a}_{d-i}$.
  By Proposition~\ref{prop:degree-quiver}, such a path exists.
  As all paths between a fixed pair of vertices of $\widetilde{Q}_\psi$ are identified in $\widetilde{A}$, the path is unique.

  The chain $p_0< \cdots < p_d$ defines a simplex in $\Order(\Path_{\widetilde{A}})$
  and so its equivalence class defines a cell in $X_{\widetilde{A}}$.
  This defines a map of sets from the dual complex to $X_{\widetilde{A}}$. Call this map $\iota$.

  To see that $\iota$ is a map of simplicial complexes, observe that the facets of the equivalence class of $p_0 < \cdots < p_d$ are exactly the simplices corresponding to chains $p_0< \cdots < \widehat{p}_j < \cdots < p_d$ for $j=0,\ldots,d$.
  For $j=1,\ldots,d$ this is exactly the image of $\{\bm{a}_0,\ldots,\widehat{\bm{a}}_{d-j},\ldots,\bm{a}_d\}$
  under the map $\iota$.
  For $j=0$ we need to show that the chain is the image of $\{\bm{a}_0,\ldots,\bm{a}_{d-1}\}$.
  The chain of paths $p_1< \cdots < p_d$ is equivalent to the chain of paths $p_1/p_1< \cdots < p_d/p_1$ in $X_{\widetilde{A}}$.
  Since each path $p_i/p_1$ is the unique path in $\widetilde{A}$ from $\bm{a}_{d-1}$ to $\bm{a}_{d-i}$, the simplex corresponding to the chain $p_1/p_1 < \cdots < p_d/p_1$ is the image of $\{\bm{a}_0,\ldots,\bm{a}_{d-1}\}$ under the map $\iota$.

  Moreover, if any pair of chains of paths $p_0< \cdots < p_d$ and $q_0< \cdots < q_d$ in $\Order(\Path_{\widetilde{A}})$ are equivalent in $X_{\widetilde{A}}$ then it must be that $h(p_i)=h(q_i)$. As $p_i$ is the unique path from $\bm{a}_d$ to $\bm{a}_{d-i}$, we have $h(p_i)=\bm{a}_{d-i}$. Thus each simplex in the dual complex maps to a unique simplex in $X_{\widetilde{A}}$ and so the map $\iota$ is injective.

  Finally, observe that the map $\sigma \mapsto e_{\iota(\sigma)}\otimes (f_{h(\iota(\sigma))}\boxtimes f_{t(\iota(\sigma))}')$ sends the generator of $S\sigma$ to the generator of $\widetilde{P}_{\iota(\sigma)}\otimes_{\widetilde{A}^e} (\widetilde{T}\boxtimes_{\field} \widetilde{T}')$ and as $\iota$ is a map of simplicial complexes, this map will be compatible with the boundary maps and so gives a well-defined map of complexes. Thus the map is an inclusion of complexes.
\end{proof}

With this setup, we can move to the proof of Theorem~\ref{thm:hpa-pdim}.

\begin{proof}[Proof of Theorem~\ref{thm:hpa-pdim}]
  By Lemma~\ref{lem:inclusion}, the complex $F_{\psi^{dual}}$ includes into the complex $\widetilde{F}^{FH}\otimes_{\widetilde{A}^e} (\widetilde{T}\boxtimes_{\field} \widetilde{T}')$.
   Lemma~\ref{lem:FHquotient} then implies that $\pi_{\Lambda(L)}(F_{\psi^{dual}})$ a subcomplex of $\FFH\otimes_{A^e} (T\boxtimes_{\field} T')$.

  For the homology statement, as $\widetilde{F}^{FH}$ is a projective resolution of $\widetilde{A}$,
  \[H_0(\widetilde{F}^{FH}\otimes_{\widetilde{A}^e} (\widetilde{T}\boxtimes_{\field} \widetilde{T}')) = \Tor_0(\widetilde{A},\widetilde{T}\boxtimes_{\field} \widetilde{T}') = \widetilde{A}\otimes_{\widetilde{A}^e} (\widetilde{T}\boxtimes_{\field} \widetilde{T}').\]

  Now define the $S$-module
  \[M := \widetilde{A}\otimes_{\widetilde{A}^e} (\widetilde{T}\boxtimes_{\field} \widetilde{T}') = \widetilde{T}'\otimes_{\widetilde{A}} \widetilde{A} \otimes_{\widetilde{A}} \widetilde{T}\cong \widetilde{T}'\otimes_{\widetilde{A}} \widetilde{T}.\]
  Then it suffices to show that $M\cong M_{\psi^{BC}}$.

  First we show that $M$ is generated as an $S$-module by pure tensors of the form $f_{\bm{v}}'\otimes f_{\bm{v}}$ where $\bm{v}$ is a vertex of the quiver $\widetilde{Q}_\psi$.
  Let $\bm{v}$ and $\bm{w}$ be distinct vertices of $\widetilde{Q}_\psi$ and let $\id_{\bm{w}}\in \End(\widetilde{T})_{\bm{0}}\cong \widetilde{A}$
  be the map on $\widetilde{T}$ which is the identity on the summand corresponding to $\bm{w}$ and zero elsewhere. Then
  \[f_{\bm{v}}'\otimes f_{\bm{w}}=f_{\bm{v}}'\otimes \id_{\bm{w}} f_{\bm{w}} = f_{\bm{v}}'\circ \id_{\bm{w}} \otimes f_{\bm{w}} = 0\otimes f_{\bm{w}} = 0.\]
  Consequently, the pure tensors corresponding to pairs of distinct vertices vanish. Thus $M$ is generated as an $S$-module by the pure tensors $f_{\bm{v}}'\otimes f_{\bm{v}}$.

  Define an $S$-linear map $M\rightarrow \field[\ZZ^{2n}]$ on the aforementioned generating set of pure tensors by $f_{\bm{v}}'\otimes f_{\bm{v}} \mapsto \frac{x^{\bm{v}}}{y^{\bm{v}}}$.
  This map is well-defined from the universal property of tensors. As the generators of $M$ map to the generators of $M_{\psi^{BC}}$,
  it follows that this map surjects onto $M_{\psi^{BC}}$. To see that this map is injective, it suffices to check on $\ZZ^{2n}$ graded components.

  Fix a degree $(\bm{u},\bm{u}')\in \ZZ^{2n}$; it suffices to show that the degree $(\bm{u},\bm{u}')$ component of $M$ is of dimension at most $1$ over $\field$.
  For a vertex $\bm{v}$, the generator $f_{\bm{v}}'\otimes f_{\bm{v}}$ is in degree $(\bm{v},-\bm{v})$ and thus for it to contribute to the degree $(\bm{u},\bm{u}')$ component requires $-\bm{u}'\leq \bm{v}\leq \bm{u}$.
  Thus to show that the component is of dimension at most $1$, it suffices to show that
  $(x^{\bm{u}-\bm{v}}f_{\bm{v}}')\otimes (f_{\bm{v}} x^{\bm{v}+\bm{u}'})=(x^{\bm{u}-\bm{w}}f_{\bm{w}}')\otimes (f_{\bm{w}} x^{\bm{w}+\bm{u}'})$ for any pair of vertices $\bm{v},\bm{w}$ with $-\bm{u}'\leq \bm{v},\bm{w}\leq \bm{u}$.

  We start with the case where there is a directed path from $\bm{v}$ to $\bm{w}$. Suppose $\bm{v}$ and $\bm{w}$ are vertices of $\widetilde{Q}_\psi$ with $-\bm{u}'\leq \bm{v},\bm{w}\leq \bm{u}$, such that there exists a path $\bm{v}\xrightarrow{\gamma} \bm{w}$. As an abuse of notation, identify the path with the element of $\End(\widetilde{T})_{\bm{0}}\cong \widetilde{A}$. Then we get the desired identity:
  \begin{align*}
    (x^{\bm{u}-\bm{v}}f_{\bm{v}}')\otimes (f_{\bm{v}} x^{\bm{v}+\bm{u}'}) &= (x^{\bm{u}-\bm{w}}f_{\bm{w}}')\circ \gamma \otimes (f_{\bm{v}} x^{\bm{v}+\bm{u}'})\\
    &= (x^{\bm{u}-\bm{w}}f_{\bm{w}}')\otimes \gamma (f_{\bm{v}} x^{\bm{v}+\bm{u'}})\\
    &= (x^{\bm{u}-\bm{w}}f_{\bm{w}}')\otimes (f_{\bm{w}} x^{\bm{w}+\bm{u}'}).
  \end{align*}

  Finally, as $(\Lambda(L)_{\RR})_{\leq (\bm{u},\bm{u}')}$ is contractible and thus connected, for any pair of vertices $\bm{v}$ and $\bm{w}$ satisfying $-\bm{u}'\leq \bm{v},\bm{w}\leq \bm{u}$, there is an undirected path from $\bm{v}$ to $\bm{w}$. Decomposing the undirected path into a sequence of directed paths reduces the question to the directed case. Thus the component of $M$ in degree $(\bm{u},\bm{u}')$ has dimension at most $1$. Thus the map is injective and so gives the desired isomorphism $M\cong M_{\psi^{BC}}$.
\end{proof}

As a final remark, we note that the constructions here do not require the conditions on the collection of line bundles that the construction of \cite{FH25}*{Section~5} requires. In particular, the line bundles in $T^{\sim}$ on $X$ will not be a full strong exceptional collection except in the simplest cases.

\section{Acknowledgments}

This paper would not be possible without the AIM workshops on ``Homological mirror symmetry and multigraded commutative algebra'' in 2025 and ``Syzygies and mirror symmetry'' in 2023. In particular, the presentation by David Favero on the results and ideas from \cite{FH25}.
The author would also like to thank Reginald Anderson, Christine Berkesch, Michael Brown, Lauren Cranton Heller, Jesse Huang, and Mahrud Sayrafi for various helpful conversations during the conception and drafting of this paper as well as feedback on earlier versions of this paper.

\subsection*{AI Usage Statement}

This paper incorporates changes and fixes based on feedback, primarily on spelling and grammar, from Claude Opus and Sonnet versions 5.5 and 5. None of the text of this paper was written by any AI tool and none of the mathematical content is sourced from any AI tool. The author takes full responsibility for all content and mathematical correctness.


\end{document}